\documentclass[a4paper,fleqn]{cas-dc}

\usepackage[numbers]{natbib}
\usepackage{soul}
\usepackage{tikz}
\usepackage{tikz,pgfplots}
\usetikzlibrary{decorations.pathreplacing}
\pgfplotsset{compat=newest}
\usepgfplotslibrary{colorbrewer}
\usetikzlibrary{arrows.meta,
                chains,
                positioning,
                shapes.geometric
                }
 
\usepackage{colortbl}
\usepackage{pgfplots}
\usepackage{pgfplotstable}
\usepackage{float}
\usepackage{array}
\usepackage{varwidth}
\usepackage{subcaption}
\usepackage{url}

\definecolor{myLightGray}{RGB}{191,191,191}
\definecolor{myGray}{RGB}{160,160,160}
\definecolor{myDarkGray}{RGB}{144,144,144}
\definecolor{myDarkRed}{RGB}{167,114,115}
\definecolor{myRed}{RGB}{255,58,70}
\definecolor{myRed2}{RGB}{215,25,28}
\definecolor{myGreen}{RGB}{0,255,71}
\definecolor{myBlue}{RGB}{44,123,182}

\newcommand{\srs}[1]{\textcolor{black}{#1}}
\newcommand{\sr}[1]{\textcolor{black}{#1}}

\newcommand{\psz}[1]{\textcolor{black}{#1}}
\newcommand{\hr}[1]{\textcolor{black}{#1}}
\newcommand{\hrblank}[1]{\textcolor{white}{#1}}

\def\tsc#1{\csdef{#1}{\textsc{\lowercase{#1}}\xspace}}
\tsc{WGM}
\tsc{QE}
\tsc{EP}
\tsc{PMS}
\tsc{BEC}
\tsc{DE}

\begin{document}
\let\WriteBookmarks\relax
\def\floatpagepagefraction{1}
\def\textpagefraction{.001}
\shorttitle{}
\shortauthors{HO Riddervold et~al.}

\title [mode = title]{A gradient boosting approach for optimal selection of bidding strategies in reservoir hydro}                      

\tnotetext[1]{}

\tnotetext[2]{}

\author[label1,label2]{Hans Ole Riddervold}[type=editor,
            auid=000,bioid=1,
            prefix=,
            role=,
            orcid=]
            
\address[label1]{NTNU}
\address[label2]{Norsk Hydro\fnref{label4}}

\cormark[1]

\ead{hans.o.riddervold@ntnu.no}
\ead[url]{https://www.ntnu.no/ansatte/hans.o.riddervold}
\credit{Idea, Writing, Original draft preparation,  Model development, Methodology}

\author[label5]{Signe Riemer-Sørensen}
\address[label5]{SINTEF}
\ead{signe.riemer-sorensen@sintef.no}
\credit{Machine learning expertise, Model development, Methodology, Software, hyper-paramater tuning}

\author[label2]{Peter Szederjesi}
\ead{peter.szederjesi@hydro.com}
\credit{Feature engineering}

\author[label1]{Magnus Korpås}
\ead{magnus.korpas@ntnu.no}
\credit{Supervision, power market theory and design }

\cortext[cor1]{Corresponding author}


\begin{abstract}
Power producers use a wide range of decision support systems to manage and plan for sales in the day-ahead electricity market, and \sr{they are often faced with the challenge of choosing the most advantageous} bidding strategy for any given day. The optimal solution is not known \srs{until} after spot clearing.  Results from the models and strategy used, and \srs{their} impact on profitability, can either continuously be registered, or simulate\srs{d} with use of historic data.  Access to an increasing amount of data opens \srs{for the application of} machine learning models \srs{to} predict the best combination of models and strategy for any given day.
In this article, historic\psz{al} performance of two given bidding strategies over several years have been analyzed with a combination of domain knowledge and machine learning techniques \srs{(gradient boosting and neural networks)}. A wide range of \hr{variables} accessible to the models prior to bidding have been evaluated to predict the optimal strategy for a given day. Results indicate that a \srs{machine learning model can learn to slightly} outperform a static strategy where one bidding method is chosen based on overall historic performance. 

\end{abstract}



\begin{keywords}
reservoir hydro \sep bidding st\srs{ra}tegies \sep hydro power \sep gradient boosting \sep neural network 
\end{keywords}

\maketitle

\section{Introduction} \label{sec1}

One of the main tasks for an operator of hydro electric power in a deregulated market is to decide how much power should be produced the \srs{following} day. Several strategies for bidding available production exist and have been described in the literature \cite{AASGARD2016181,YAMIN2004227}. 
Each strategy will potential\psz{ly} lead to different commitments for production, which again will have \psz{an} impact on profitability.

Market prices and inflow for the next day \srs{are} uncertain, and the profit associated with any selected strategy is not revealed \srs{until} after the decisions are made. 

The question address\srs{ed} in this article, is if the producer with access to sufficient amounts of information about historic\psz{al} performance of different strategies, can predict in advance which bidding strategy should be selected for a given day. 

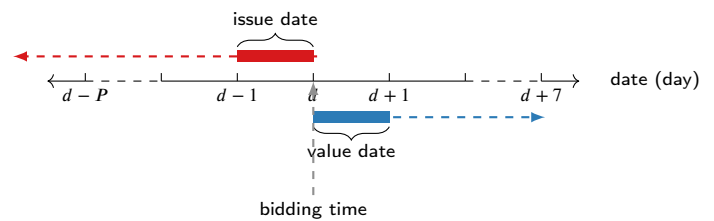
\begin{figure}[pos=htbp]
\begin{center}

\begin{tikzpicture}[%
    every node/.style={
        font=\scriptsize,
        text height=1ex,
        text depth=.25ex,
    },
]
\draw[<-] (1.5,0) -- (2,0);
\draw[dashed] (2,0) -- (3,0);
\draw[-] (3,0) -- (7,0);
\draw[dashed] (7,0) -- (8,0);
\draw[->] (8,0) -- (8.5,0);

\foreach \x in {2,3,...,8}{
    \draw (\x cm,3pt) -- (\x cm,0pt);
}

\node[anchor=north] at (2,0) {$d-P$};
\node[anchor=north] at (4,0) {$d-1$};
\node[anchor=north] at (5,0) {$d$};
\node[anchor=north] at (6,0) {$d+1$};
\node[anchor=north] at (8,0) {$d+7$};
\node[anchor=north] at (9.5,0.2) {date (day)};

\fill[myRed2] (4,0.25) rectangle (5,0.4);
\draw[myRed2,dashed,thick,-latex] (5.05,0.325) -- (1.05,0.325);

\fill[myBlue] (5,-0.4) rectangle (6,-0.55);
\draw[myBlue,dashed,thick,-latex] (6.05,-0.475) -- (8.05,-0.475);

\draw[decorate,decoration={brace,amplitude=5pt}] (4,0.45) -- (5,0.45)
    node[anchor=south,midway,above=4pt] {issue date};
\draw[decorate,decoration={brace,amplitude=5pt}] (6,-0.6) -- (5,-0.6)
    node[anchor=north,midway,below=4pt] {value date};
\draw[myDarkGray,dashed,thick,-latex] (5,-1.5) -- (5,0);
\node[anchor=north] at (5,-1.5) {bidding time};

\end{tikzpicture}
\end{center}

\caption{Definition of issue- and value date where ``d'' is date, and ``P'' is any date before the bidding time.} \label{fig:issuevaluedate}
\end{figure}

Bidding \hr{day-ahead production} to the power-exchange is typically done the day before the actual commitments are executed. 
"Issue date" is defined \sr{as} the date when bidding is conducted, while "value date" is used for the date when commitments from the bidding \psz{are} realized through costs and income. 
Only variables that can be identified on or before the issue date, can be used to classify the optimal strategy associated with performance of a corresponding value date. Fig. \ref{fig:issuevaluedate} describes the relationship between the two.

\hr{In the following sections, results from a case study  investigating historical performance of two different bidding strategies will be presented \cite{Riddervold8916227}. Further, the steps associated with the machine learning process applied in this article will be described together with examples of application on the historical data. In Section \ref{casestudy}, a concept for application and case-study will be presented, before drawing a conclusion in Section \ref{Conclusion}.}

\section{Results from historic\psz{al} bidding strategies} \label{historic}

Two strategies have been evaluated in this analysis. \sr{In the first method, the} expected volumes are found by deterministic optimization with price forecast and inflow 
and 
submitted as fixed hourly bids to the power exchange. \sr{The optimization is performed with SHOP, which} is a software tool for optimal short-term hydropower scheduling developed by SINTEF Energy Research, and is used by many hydropower producers in the Nordic market \cite{FossoO.B2004Shsi}. The second strategy is stochastic bidding. The stochastic model is based on the deterministic method, but allows for a stochastic representation of inflow to the reservoir and \sr{the} day-ahead market prices. In this case, bid-curves can be generated from \psz{the} stochastic model as described in \cite{AasgardEllen2018Odbc}. Example of results from evaluation of two strategies for some selected days \sr{for a specific river system} are shown in Table \ref{table:1}. The \srs{performance quantification}  for \srs{the} deterministic and stochastic \srs{models} are the performance-gaps in EUR for these strategies.

\begin{table*}[width=1.8\linewidth,cols=7,pos=h] 
\caption{Example of performance quantification \sr{in EUR} for the deterministic and stochastic models for \sr{three} days in 2017 \sr{for the use-case river system}.}\label{table:1}
\scalebox{0.9}{%
\begin{tabular}{lc c c c c c c}
\toprule
\sr{I}ssue date  & \sr{V}alue date & Deterministic ($\beta_{det\psz{}}$) 
& Stochastic ($\beta_{stoch}$) 
& MIN	& DELTA($\eta_{s}$) & BEST\\
\midrule
2017-07-01 & 2017-07-02 & 69.2	& 137.9	& 69.2	& 68.7 & 1 \\
2017-07-02 & 2017-07-03	& 16.5 & 65.1	& 16.5	& 48.6	& 1	\\
2017-07-03 & 2017-07-04 & 31.1 & 29.9	& 29.9 & -1.2 & 0 \\
\bottomrule
\end{tabular}}
\end{table*}

In \cite{Riddervold8916227}, a method for 
\hr{measuring performance of individual historic\psz{al}} bidding days has been proposed, where \hr{$\Pi_{s,d}$ in Eq.\ref{eq:1} represents a measure of value for bidding strategy {\it s} on day {\it d}, and $\Pi_{opt,d}$ describe the optimal value for the relevant bidding date based on a deterministic strategy with perfect foresight of price and inflow}. 



\sr{The p}erformance gap ($\beta_{s,d}$) reflect\sr{s} the loss of choosing strategy \textit{s} for day \textit{d} compared to an optimal deterministic strategy, and is calculated as the difference between $\Pi_{opt,d}$ and $\Pi_{s,d}$. A high \sr{value} for $\beta_{s,d}$ indicate poor performance.

\begin{equation}\label{eq:1}
\beta_{s,d} = \Pi_{opt,d}-\Pi_{s,d}    
\end{equation}
\begin{equation}\label{eq:2}
\eta_{s,d} =  \beta_{stoch,d} - \beta_{det,d}  
\end{equation}

\sr{Since the deterministic and stochastic predictions are based on pre-bidding values, even the best model will rarely correspond exactly to the perfect foresight strategy computed after the actual inflows and prices are known. Consequently, we define the best model as the one with the lowest gap relative to the perfect scenario. We classify each date (data point) with ``0'' if the stochastic strategy leads to to lowest gap, and ``1'' if the deterministic strategy leads to the lowest gap}
\hr{as shown in table \ref{table:1}}. 

If we define ``strategy gap'' ($\eta_{\hr{s,d}}$) as ``performance gap stochastic'' ($\beta_{stoch,d}$) minus ``performance gap deterministic'' ($\beta_{det\psz{},d}$), a high value 
\sr{is} a strong signal to choose a deterministic model for that day. Negative values indicate that a stochastic model is preferred, and the more negative, the higher 
the importance of a \hr{stochastic} 
model. Values around zero indicate \sr{a negligible difference in the choice of model for the day.} 
With this insight, we see that a pure classification model only measuring the correct number of classifications, will not necessarily give the best results if the overall target is to have a low deviation from optimum over time.
 
\sr{It seems obvious to apply supervised machine learning to the selection problem.}
\sr{We have tested two different approaches: i) labelling each data point as stochastic or deterministic based on the performance gap, and then training a classification model to predict which category unseen data points belong to. 2) A regression model trained to predict the performances of each model directly and using a simple decision heuristic (minimum gap) to decide on the most advantageous strategy.}

In Sec.\ref{ML}, we will systematically go through the different steps associated with the machine learning challenges for the two approaches.

\section{Machine Learning Process} \label{ML}
The use of machine learning to classify strategies or predict values have gained significant momentum during the last few years \cite{Jordan255}. Machine learning is a set of techniques that allow a 
\sr{a computer algorithm} to improve performance as it gains experience, which in our case is exposure to additional data. No explicit instructions are required, but instead, the algorithm needs some sort of training on representative input and output data. 
If the training is successful and the model can find some general patterns or behaviour in the data, it can subsequently be used to \sr{predict output} for previously unseen input data. Machine learning basically covers everything from simple regression to deep neural networks. 

Classification and regression problems, are 
categorised as supervised learning methods where the training takes place on pairs of input and output data, and subsequently applied to unseen input data to generate predictions. Other categories are unsupervised learning used for clustering and association, and reinforcement learning used for decision optimisation. 

Within the area of electricity power market analysis, neural networks have been used to investigate strategic and algorithmic bidding \cite{8424030,PINTO201627} as well as for price- and load forecasting \cite{ANBAZHAGAN2012140,LI201596}. 

Gradient boosting methods have received less attention, but have successfully been applied for load- and price forecasting \cite{7885697,IRIA2019324}.

Results from the previously published articles, indicated that machine learning techniques successfully can be applied to improve forecasting of load and prices in the power market, but there exists limited publications documenting operational use and added value from these techniques compared to what is state of the art in the industry. Several energy companies as well as software-, data- and consultancy companies supplying the energy sector advertise and promote the use of machine learning, and the European intra-day market has in particular been an area of interest in relation to application \cite{Wattsight,engie}. 

\sr{In this article, we assume that the data is time-independent in the sense that a strategy choice for the next day does not affect which strategy will be the best choice further in the future. Such time dependence could potentially be accounted for in a reinforcement learning framework which has received increased attention as a solution to the increasing amount of non-linear relationships and high-dimension problems associated with hydropower production planning \cite{10.1007/978-3-030-22999-3_11, 10.1007/978-3-319-71078-5_15, Zarghami2018ShortTM}. However, in general, reinforcement methods are still fairly immature and require significant fine tuning on individual problems in order to work \cite{rlblogpost}, so we consider it outside the scope of this work.}

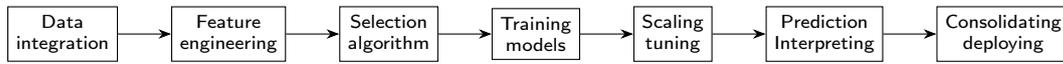
\begin{figure*}
\centering
 \begin{tikzpicture}[
    node distance = 5mm and 7mm,
      start chain = going right,
 disc/.style = {shape=cylinder, draw, shape aspect=0.3,
                shape border rotate=90,
                text width=12mm, align=center, font=\linespread{0.5}\selectfont},
  mdl/.style = {shape=ellipse, aspect=2.2, draw},
  alg/.style = {draw, align=center, font=\scriptsize}
                    ]
    \begin{scope}[every node/.append style={on chain, join=by -Stealth}]
\node (n1) [alg] {Data \\ integration};
\node (n2) [alg]  {Feature \\ engineering};
\node (n3) [alg]  {Selection\\ algorithm};
\node (n4) [alg] {Training \\ models};
\node (n5) [alg]  {Scaling \\ tuning};
\node (n6) [alg]  {Prediction \\ Interpreting};
\node (n7) [alg]  {Consolidating  \\ deploying};

    \end{scope}
    \end{tikzpicture}
\caption{Typical sub-tasks involved in machine learning analysis process.} \label{fig:MLprocess}
\end{figure*}

There are several ways of approaching a machine learning problem, but they often involve some or all of the sub-tasks illustrated in Fig. \ref{fig:MLprocess} 
and several iterative loops over the process. Our implementation is discussed in Sec. \ref{sec:DataIntegration}-\ref{sec:Prediction}.


\subsection{Data integration} \label{sec:DataIntegration}

\subsubsection{Data gathering} \label{ssec:Datagathering}
Data gathering can be time consuming and the data might require significant cleaning and quality assurance before it can be injected into a learning model. In addition to gathering data for the values we want to predict, defined as output variables, it could be wise to have an idea about what could be relevant input variables to be collected in the same process. In the case of predicting the best strategy using a classification model, the output variable is the prediction of a stochastic or deterministic model, while the input variables could be inflow and prices. The output from a regression model would be the strategy gap ($\eta_{s}$). 
That is, given some input variables (input), what is the predicted output variable (output). 


Experienced production planners \sr{may} have an \sr{intuitive} 
perception of what the important input variables could be, and these can be used as initial values in the learning process. \sr{This is often} an iterative process where additional insight gradually is gained and \sr{frequent reviews of the initial conditions are required. In order to maximise learning and minimise bias, it can be necessary to adjust the variables.} Insufficient performance of \sr{a} model \sr{should call for}
rethinking both related to input, model \srs{structure} and system-design.

\subsubsection{Input variables} 
The two strategies evaluated in this analysis are measured against a perfect foresight model where prices and inflow are known prior to bidding. All other factors going into the performance evaluation are equal. We expect that factors affecting the prices and inflows will have a significant effect on the performance gap. 
Another important factor affecting the production schedule is the water value. The water value can be defined as the future expected value of the stored marginal kWh of water, i.e. its alternative cost \cite{WOLFGANG20091642, BROVOLD2014117}.

The basic 
\hr{concept} is to produce \sr{when} the water value is lower than the price. The lower the water value is compared to the price, the 
\hr{stronger is the signal to produce}. 
Based on this insight and domain knowledge of what might effect production, 
\sr{we apply some} general hypotheses 
to select an initial set of input variables \sr{for further analysis as given in Table \ref{table:2}.} 







\begin{table*}[width=1.9\linewidth,cols=5,pos=h]
\caption{Initial set of input variables, hypotheses for how they affect the choice of strategy, and the strategy that the \sr{domain} experts presumably would select for high/low values of the variables.} \label{table:2}
\scalebox{0.8}{%
\begin{tabular}{lp{3cm}p{13cm}p{0.4cm}p{0.4cm}}

\toprule
Nr. & Variable & Hypothesis & High & Low   \\
\midrule
1 & Inflow deviation from \hr{historic} normal & Higher observed inflow increase\psz{s the} risk of flooding for the next day (value date) & S & D  \\
2 & Reservoir filling & High or low reservoir filling increase\psz{s the} risk of flooding or resource shortage during value date & S & S  \\
3 & Price volatility & High volatility in prices gives increased uncertainty for prices the next day & S & D \\
4 & Price volatility in \psz{the} prognosis & High volatility in \psz{the} price prognosis give increased uncertainty for prices the next day & S & D \\
5 & Water value  & Water value is the primary \sr{deciding factor for production,} 
and will clearly have an effect on the \sr{strategy} choice when seen in relation to other variables  & U & U  \\
6 & Average price & Price relative to water value is important  & U & U \\
7 & Average price prognosis & Price prognosis relative to water value is important & U & U\\
\bottomrule
\end{tabular}}
\end{table*}

All input variables in \srs{Table \ref{table:2} are related to the issue data, and consequently available when the strategy for the next day is decided}. 
In the high and low columns, it is indicated what is believed to be the preferred strategy based on experience from operation in case of high and low values. ``U'' indicates that the domain experts are uncertain how the value will affect the classification, but still believe it is important relative to choice of strategy.

\subsubsection{Preparation and visualisation} 
\begin{figure}[pos=htbp]
	\centering
		\includegraphics[width=0.8\columnwidth]{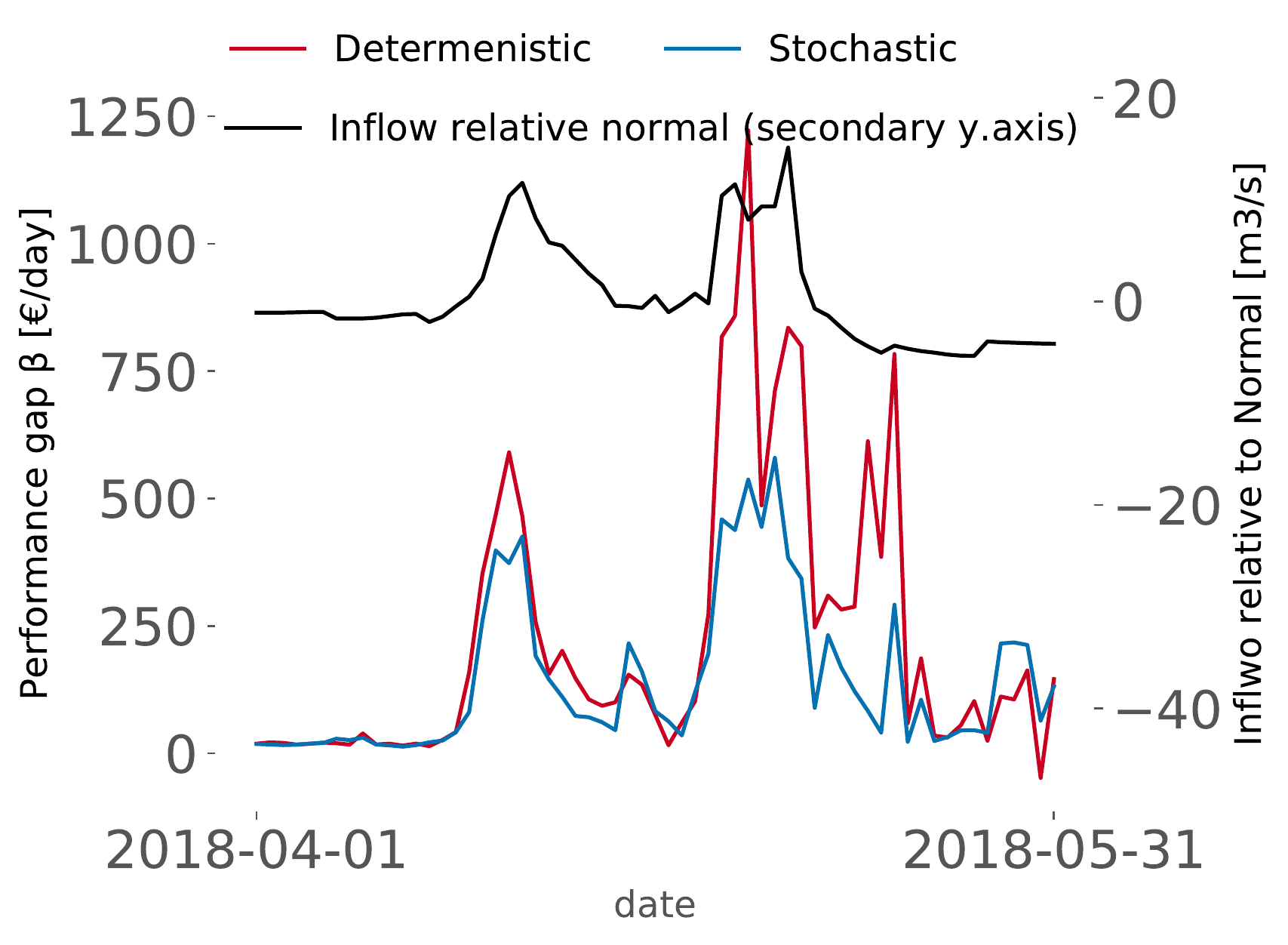}
	\caption{Detailed results for performance gap from April-May 2018, compared to inflow relative to normal. Data from \cite{Riddervold8916227}.}
	\label{fig:performacegap_vs_inflow}
\end{figure}

Visualisation and sanity check of data is important in any data analysis.
The first step is to visualize all the input- (X) and output variables (y). This is to detect missing or obvious\psz{ly} incorrect data in the data set, or anything that should be corrected for during the analysis. A plot \sr{may quickly reveal periods of missing data or unrealistic values, or provide ideas of interactions between input and output values.} 
Fig. \ref{fig:performacegap_vs_inflow} \sr{shows performance gaps ($\beta$) together with inflow deviation from normal. It provides a clear indication of a connection between periods with high inflow and poor performance of deterministic bidding.} 

\begin{figure}[pos=htbp]
    \centering
    \includegraphics[width=0.8\columnwidth]{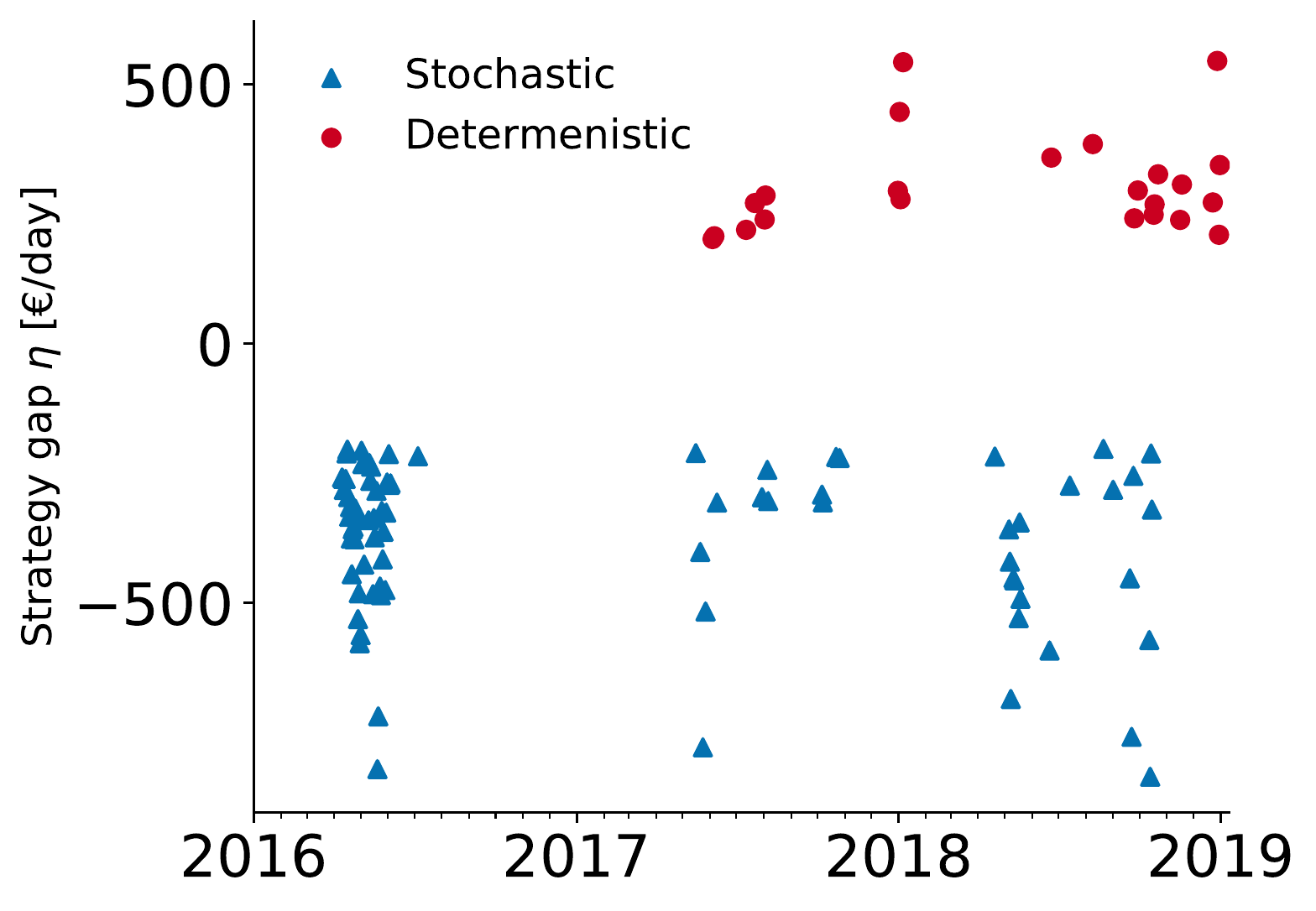}
    \caption{Strategy gaps in the period 2016-2108.} 
    \label{fig:gapsover200}
\end{figure}

 \sr{In Fig. \ref{fig:gapsover200} we plot the strategy gap ($\eta_{s}$) defined in Eq. \ref{eq:2}, but only for absolute values larger than 200 €/day corresponding to days where there is a significant impact from the choice of strategy. These values typically concentrate around the second quarter every year, indicating that time of year e.g. month may be a relevant variable to include in the analysis.} 
 

\subsubsection{Correlation} \label{ssect:Correlation} 
\begin{figure}[pos=htbp]
\centering
  \includegraphics[width=0.8\columnwidth]{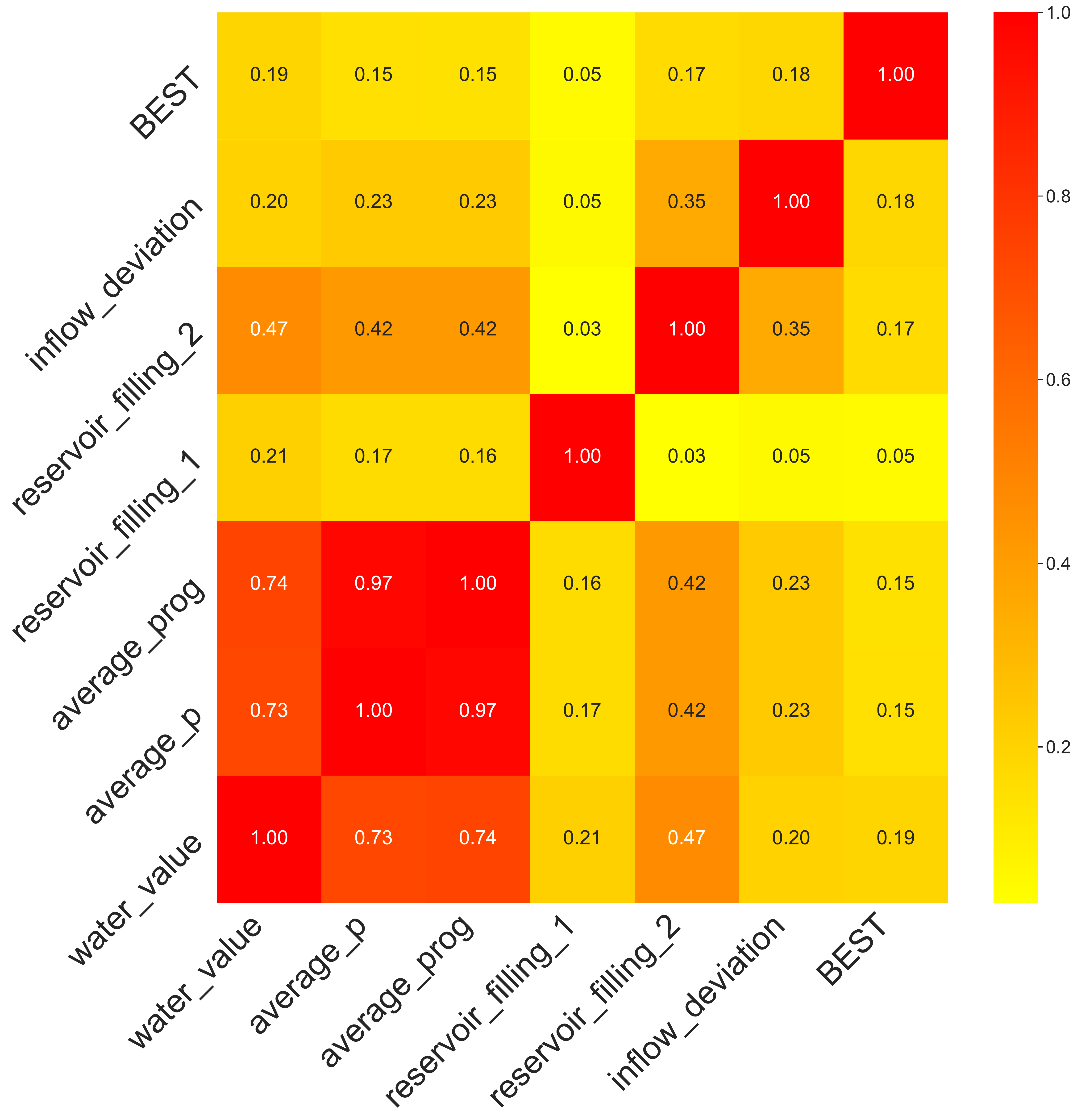}
    \caption{Correlation matrix for the initial set of variables.}
  \label{fig:matrix}
\end{figure}

In Fig. \ref{fig:matrix} \srs{relative (Spearman)} correlation between variables are shown.
If some of the independent variables (X) are un-correlated with the dependent variable(s) (y) of interest, they can be removed. If some of the independent 
variables are correlated with each other, \srs{a subset can be removed} or, even better, combined.


We are mainly interested in the correlation between variables and the prediction indicated in the ``BEST'' column. A higher number indicates that the variable will play a more important role in the classification. In this case, the water value and reservoir filling in reservoir 2 are the variables with the highest correlation to the prediction. There is also a strong correlation between prices associated with the issue date (\texttt{average\_p}) and the prognosis for the value date (\texttt{average\_prog)}.

\begin{figure}[pos=htbp]
  \centering
  \includegraphics[width=0.8\columnwidth]{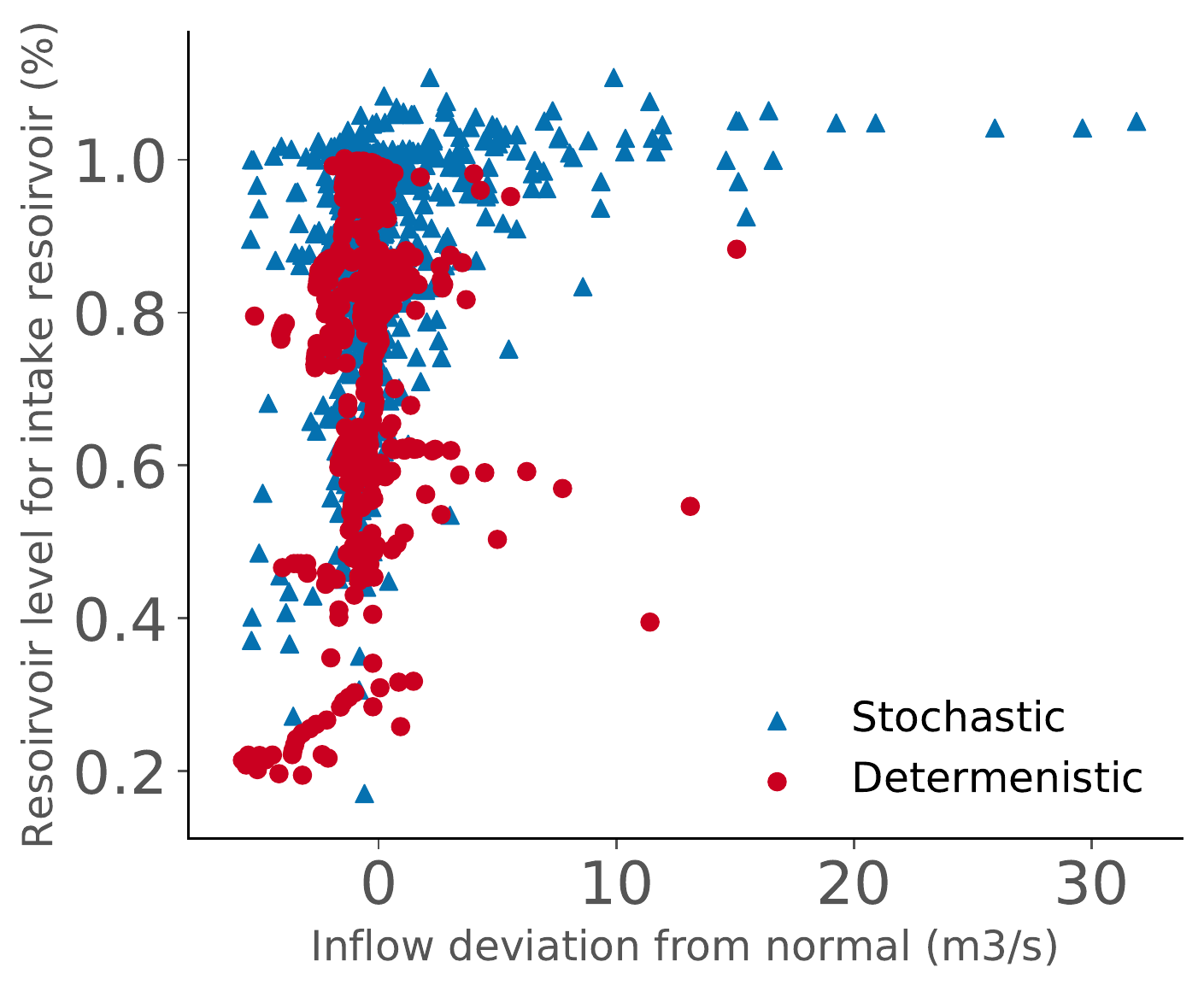}
    \caption{Gridplot for variables inflow and reservoir \sr{level (denoted reservoir filling in Fig. \ref{fig:matrix} ).}}
  \label{fig:gridplot}
\end{figure}

Grid plots \sr{as shown in Fig. \ref{fig:gridplot}} can also be used to plot all variables against \sr{each other} 
to see if they group clearly into categories. 
A trained eye will spot the correlations from the grid plot. In this plot \hr{blue} (\hr{triangular}) markers indicated when the stochastic model performs best, while \hr{red} (round) markers indicate when the deterministic model is preferred. When the level in the intake reservoir is high or inflow is well above normal, a stochastic model is preferred. 


                


\def\cca#1{\cellcolor{black!#10}\ifnum #1>5\color{white}\fi{#1}}



\subsection{Feature selection and engineering}
The terms variable and feature are often used without distinction when referring to machine learning problems \sr{unless for kernel methods for which features are not explicitly computed \cite{2685}}. 

In this article, input variables are primarily used when referring to variables collected in the data collection phase, while we use the term feature for constructed variables engineered to maximise the information available for the model. 

We thereby define feature engineering as selection and combination of relevant information we are investigating in order to optimally exploit the available characteristics in the data. \sr{In the ideal world with unlimited training data and computational power, the machine learning algorithm can be fed all the available variables and figure out the relevant features itself. However, with limited access to data, feature engineering is often an iterative process where you go back and forth between data pre-processing, correlations, feature engineering and model fitting.}

In this case study, two examples of extracted features are month and the strategy-gap associated with the issue date. In the input data, date is given, and specifying the month is therefore a trivial task. The same applies for strategy gap. 
In the initial choice of variables, we are not considering the time dimension of each variable. Given our physical knowledge of the system, we know there could be some delay between parameter changes. Consequently, it might be better to shift some parameters by a \srs{certain} amount of time to get a variable called X-t. As illustrated in Table. \ref{table:1}, the data contains the strategy gap for each day, and the new feature is then the strategy gap for value day minus one.

More advanced methods such as principal component analysis, can also be used to get ideas for feature engineering, but have not been exploited further in this analysis.

\subsubsection{Bid- and ask curves} \label{sssec:Bidaskcurves}
\begin{figure*}[width=2.0\linewidth] 
\centering
   \begin{subfigure}{0.46\linewidth} \centering
     \includegraphics[scale=0.36]{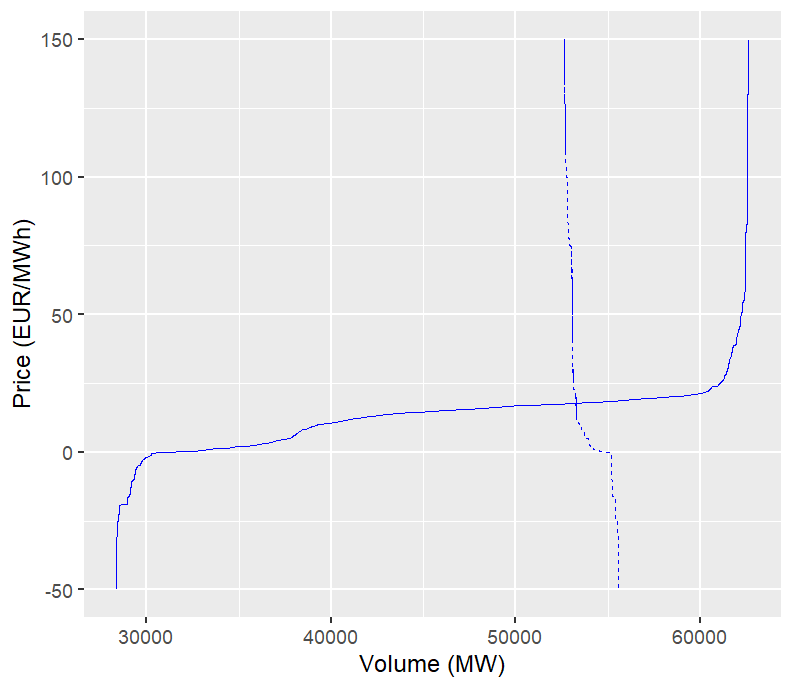}
     \caption{Bid- \sr{(dashed)} and ask \sr{(solid)} curves for February 11, 2016, hour 11}\label{fig:bidask_1}
   \end{subfigure}
   \begin{subfigure}{0.46\linewidth} \centering
     \includegraphics[scale=0.36]{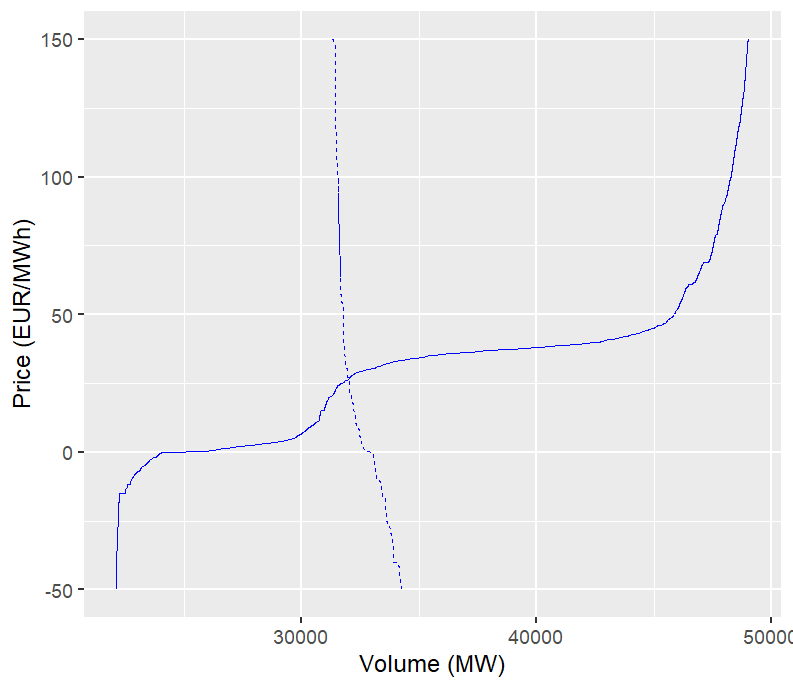}
     \caption{Bid- \sr{(dashed)} and ask \sr{(solid)} curves for May 1,  2018, hour 2 }\label{fig:bidask_2}
   \end{subfigure}
   \caption{Bid- and ask curves Nord Pool power market}
   \label{fig:bidaskdates}
\end{figure*}

\sr{It can be hypothesised} that price volatility can be predicted by investigating bid/ask sensitivities in the day-ahead market. In the Nordic Market, the spot price for each hour is determined by the crossing-point between the bid and ask curves published daily after spot clearing by Nord Pool. 
Two examples of bid- and ask curves after interpolation are illustrated in Figures \ref{fig:bidask_1} and \ref{fig:bidask_2}, where the dashed line is the bid curve and the solid line is the ask curve. In the left plot, the curves show a fairly stable balance but with a larger price sensitivity to the upside, whereas the curves in the right plot show a balance with strong sensitivity to the downside, i.e. a price collapse.  

Thus, the spot price volatility depends on the shape of the two curves and how they change from hour to hour. As seen, the \sr{(dashed)} bid curve is normally vertically shaped, implying that power consumption on a given day is fairly inelastic. The \sr{(solid)} ask curve has a significant price drop for low volumes and a significant price rise for high volumes. This is reflecting the underlying price elasticity of supply from different sources of energy.

The hypothesis is that when demand is either close to the point where the bid price increases almost exponentially, or close to the point where bid price drops toward zero, an increased volatility is probable, and a stochastic bidding strategy is preferred. Again, only information available prior to bidding can be used to predict the optimal strategy the next day, and we therefore use the bid curves for the issue date to predict strategy for the value date. 

\begin{figure*}[width=2.0\linewidth] 
\centering
   \begin{subfigure}{0.46\linewidth} \centering
     \includegraphics[scale=0.36]{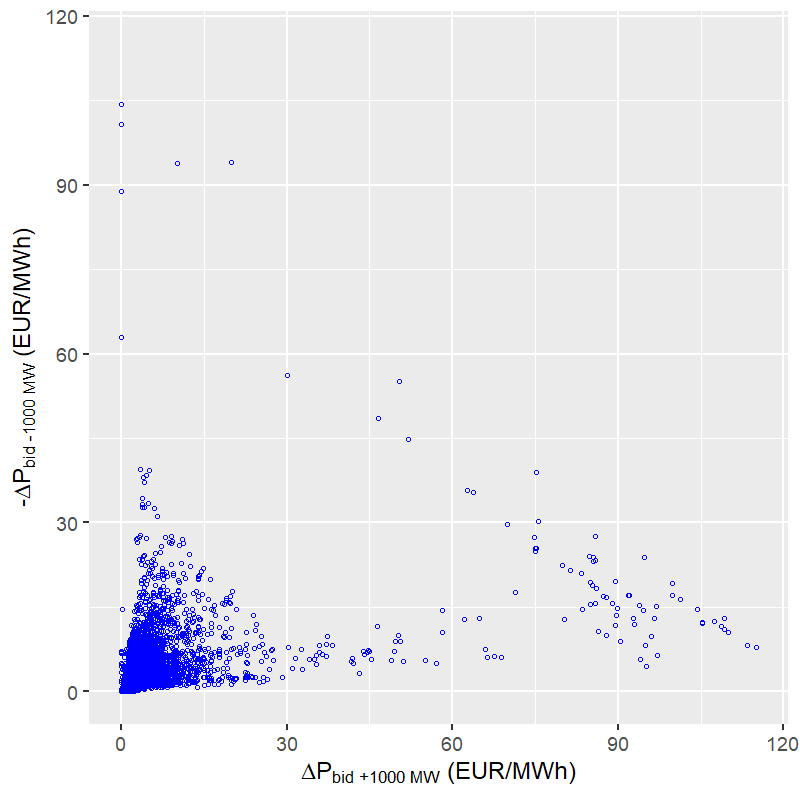}
     \caption{Price difference based on a +/- 1000 MW shift in load.}\label{fig:bidask_sensitivities}
   \end{subfigure}
   \begin{subfigure}{0.46\linewidth} \centering
     \includegraphics[scale=0.36]{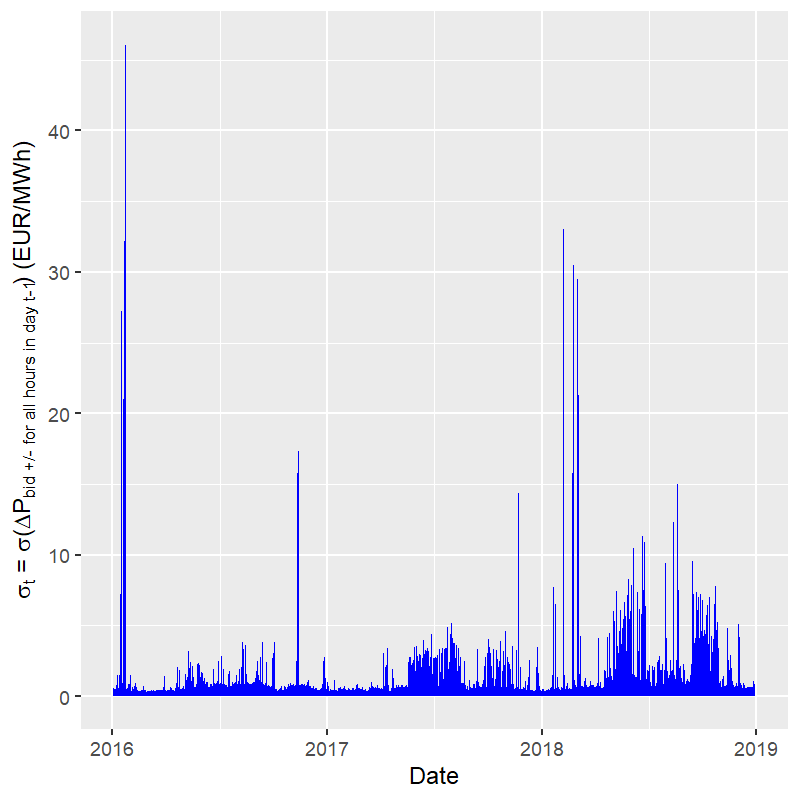}
     \caption{Rolling volatility feature based on the standard deviation of price sensitivities for the previous day.}\label{fig:bidask_volroll}
   \end{subfigure}
   \caption{Bid- and ask sensitivities and rolling volatility in the Nord Pool power market.} \label{fig:bidask}
\end{figure*}

As a base we use the bids and asks of volume for the full range of prices available by Nord Pool. The bid and ask curves are constructed by interpolating the points with a 0.1 EUR/MWh granularity. Subsequently, they are then used to compute price sensitivities based on shifts in the bid curve. The resulting features are the price differences resulting from a 1000 MW changes in load, i.e. one price difference after a horizontal shift of +1000 MW in the bid curve and one price difference based on a -1000 MW shift. The price difference is computed as the crossing-point price between the bid- and ask curves after a shift, minus the pre-shift crossing-point price. The computed hourly sensitivities are presented in Fig. \ref{fig:bidask_sensitivities}. These are used to compute a rolling volatility feature based on the standard deviation of price sensitivities for the previous day, presented over time in Fig. \ref{fig:bidask_volroll}. 


Thus, the price sensitivity to shifts in the bid curve varies over time and can occasionally experience sudden spikes or build-up over multiple days. The assumption that the next day will experience a high volatility if the current day had a high volatility is intended to be captured by the t-1 rolling volatility feature. An additional t-2 feature was tested, but did not yield any additional qualitative improvements in the classification of optimal strategy. Further, rate-of-change features were computed as the differences in standard deviation of price sensitivities from day t-2 to day t-1, but did not yield any significant improvements in the quality of predictions when added to the model.

\subsubsection{Features used}
To evaluate the benefit of increasing the amount of variables, two input data-sets are defined. The first data-set is referred to as {\it simple input}, and consists of 
eight variables \sr{listed in Table \ref{table:2}}. The second set is referred to as {\it complex input} and consists of \hr{the following variables}:
\begin{itemize}
\item \sr{the} eight \sr{simple} variables
\item all hourly prices for both issue date and prognosis for value date
\item bid-ask curves
\item rolling volatility
\item month, year, day and performance of similar week-days
\item strategy gap for issue date
\item rate of change for reservoir filling
\item difference between price and water value
\end{itemize}
In total, the complex-input data-set consist of 113 variables.

\subsection{Selecting the model and algorithm}

The problem 
\sr{at hand} can be solved as a classification problem where the target is to predict the best bidding strategy \sr{or by predicting the }
strategy gap directly, and based on this decide the optimal strategy. 

\sr{An} algorithm suitable for \sr{both approaches is boosted decision trees as implemeted in the } 
XGBoost \sr{library}  \cite{Chen2016XGBoostAS}. 
XGBoost 
implements the gradient boosting \cite{hastie_09_elements-of.statistical-learning} decision tree algorithm and use multiple trees to classify samples . \sr{This} gives valuable insight \sr{on the feature importance and possible feature interactions} 
The main hyper parameters 
are learning rate, maximum tree depth and number of estimators represented by trees \cite{XGBoost}.

\begin{figure*}[width=2.0\linewidth]
\centering
\includegraphics[width=1.8\columnwidth]{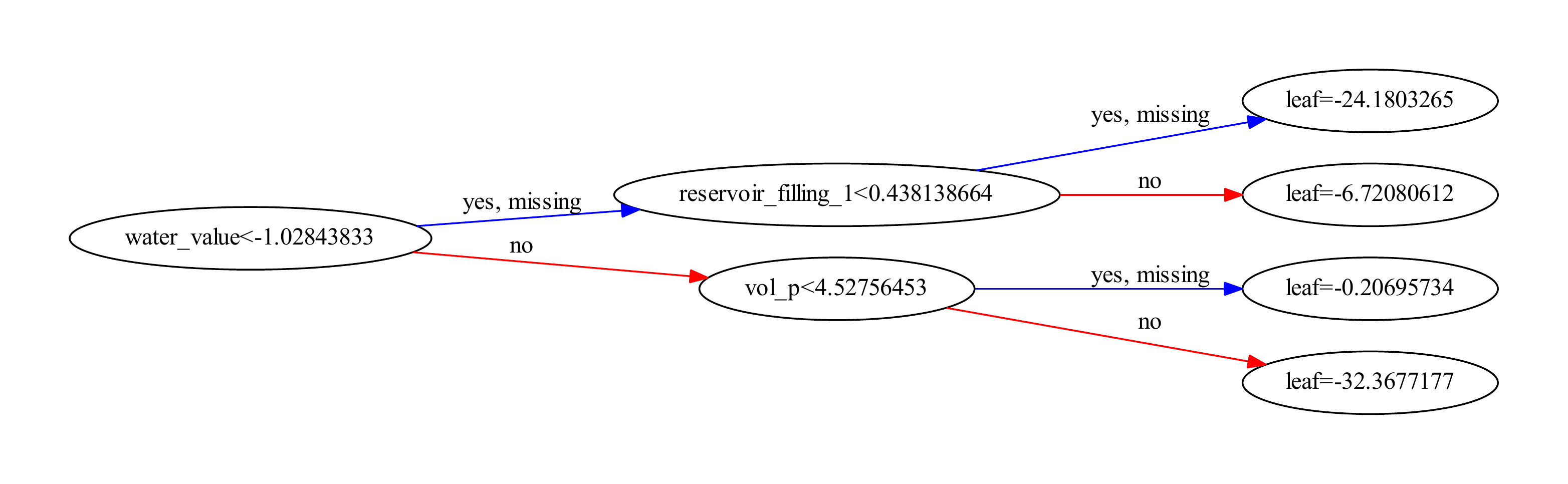}
\caption{An example of one classification tree from the model trained on a limited set of features. The total model is an average over many trees. } \label{fig:classtree}
\end{figure*}

Figure \ref{fig:classtree} illustrates a random\sr{ly selected} tree from \sr{the resulting classification model for} a limited set of features. The first split feature is placed highest or to the left in the tree, followed by subsequent criteria. The total model is an average of multiple trees with different representations of the features. 
The importance of each parameter in the trees can be derived in several ways. In this article,  we have applied two methods to quantify the importance. These are GAIN and SHAP which will be described further in Sec. \ref{ssec:featurescaling}. 




\srs{Fully connected neural networks (NN)} 
and recurrent neural networks (RNN) have been tested with less success and a comparison with these methods are included in Sec. \ref{sssec:neuralnetwork}. 


\subsection{\hr{Splitting the data}}
Before we can train a model, we need to split the data into training samples and test samples. The basic principle is that we build and train the model based on a fraction of the data called training data, and then test the performance of the model on data that are not used when building the model. \sr{This sample is} referred to as test data, and represent unseen samples that are not accessible in the model building phase. 
There exist several ways of generating training and test sets. A common approach is to define a the ratio between training and test data, e.g. 0.7/0.3, and let the model randomly split the data into sub-sets. 

Applying the concept of randomly splitting the 
data into training and test might conflict with the underlying mechanisms behind the power market. When managing Hydro power, weather plays an important role. The possible outcomes and combinations of reservoir levels, inflow, snow and prices are widely spread, and the power producers therefore tend to use a long history of observations to represent the possible scenarios in there decision support models. Typically this can be in the range of 30-60 years of history. Consequently a random split \sr{of historic data} may leak information from the \sr{historic} future into the training sample, which is not desirable. Instead, we can use a sequential split of data, and test the trained model on years not previously seen by the model. Such a split will lead to a model that performs well on years with similar historic representation, but with very little predictive power in ``freak'' years. Three years of data are available in this analysis. If we assume that two years are used for training, and one year is used for testing, we get a split of 67/33 between training and test data.  In addition to random splitting of data, splitting of data have been tested with the years 2016 and 2017 used for training, and then tested on 2018 data. Results from this analysis is presented in Sec. \ref{ssection:Classification}.





Another method which has not been tested in this analysis, is to use all available historic information, or a rolling number of historic days to train the model. The model will then be used to classify the next day. To test performance of such a model, a simulation framework would be needed since we are depending on updating the training set every day, as well as re-calibrate the model.

\subsection{Feature scaling and hyperparameter tuning}\label{ssec:featurescaling}
\subsubsection{Feature scaling}
Not all features have the same scale: Some have values of the order of 1000s, and some are 0.1. In order to let them equally influence the model, we need to ``put everything on the same scale''. We can either scale everything to a fixed range of values 
or change the distribution to become a normalised Gaussian. 

In general, decision tree algorithms such as XGBoost do not require scaling \cite{XGBoost}, but it might help with quicker convergence in relation to numerical processing. Scaling is required for neural networks, so to be able to apply the same pre-processing of data, similar scaling has been applied for both XGBoost and neural networks in this article.  

Depending on the sample size, the test data can either be scaled with their own scaling (for large samples), or with the training sample (small samples). What to chose depends on how you would treat the actual data you will later use with the model.

\begin{figure*}[width=1.9\linewidth] 
\centering
   \begin{subfigure}{0.3\linewidth} \centering
     \includegraphics[scale=0.3]{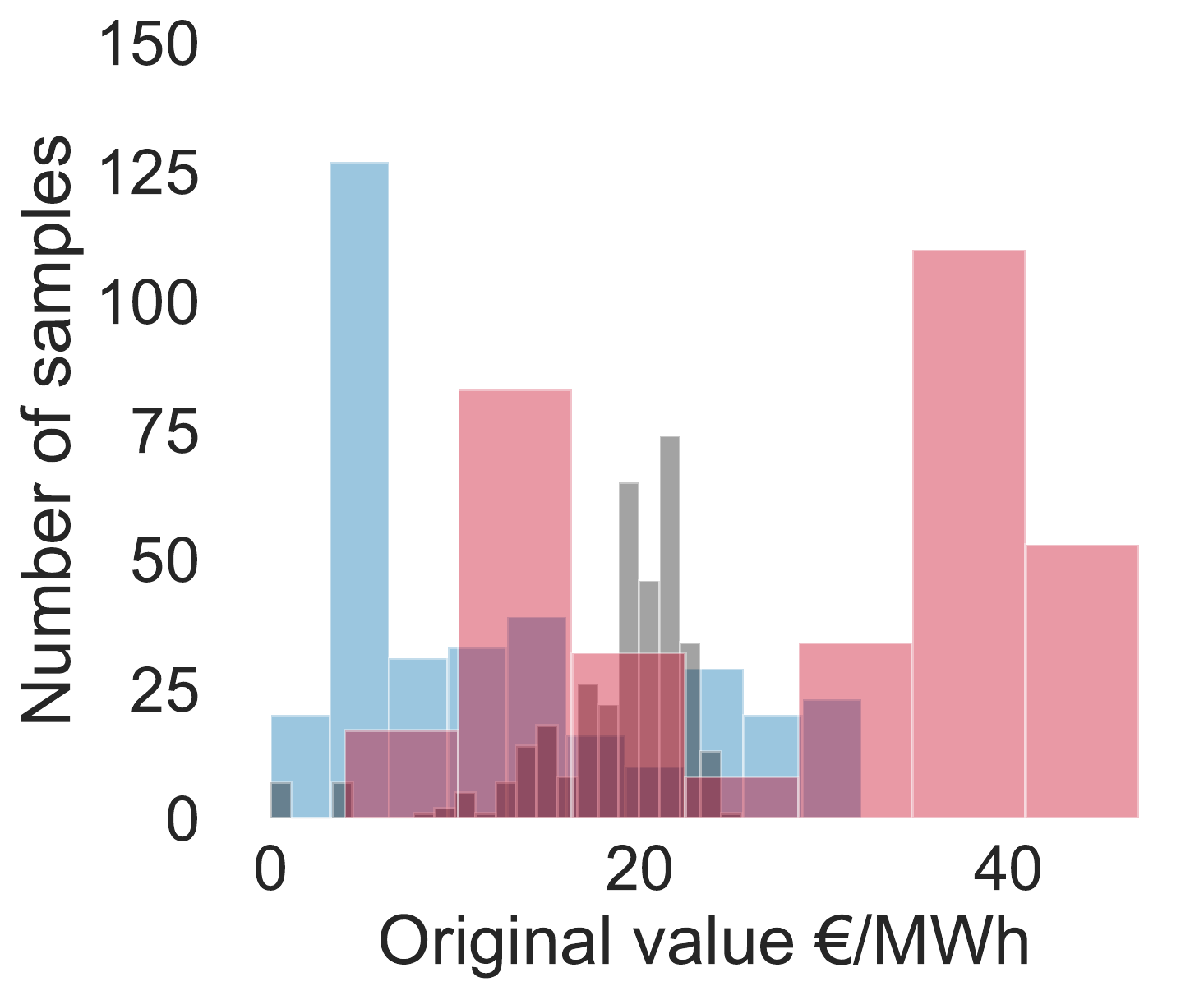}
     \caption{}\label{fig:dist1}
   \end{subfigure}
   \begin{subfigure}{0.3\linewidth} \centering
     \includegraphics[scale=0.3]{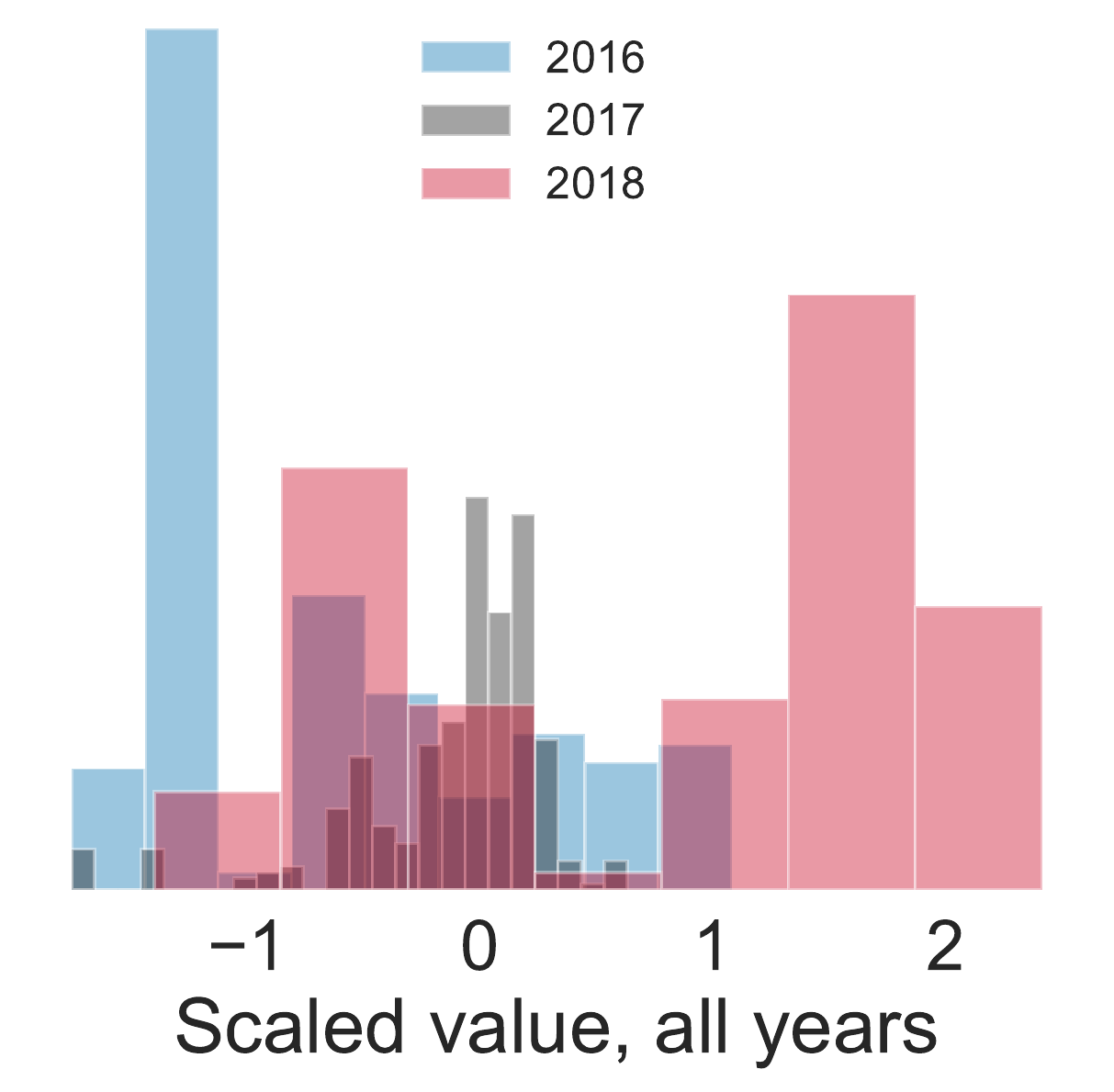}
     \caption{}\label{fig:dist2}
   \end{subfigure}
   \begin{subfigure}{0.3\linewidth} \centering
     \includegraphics[scale=0.3]{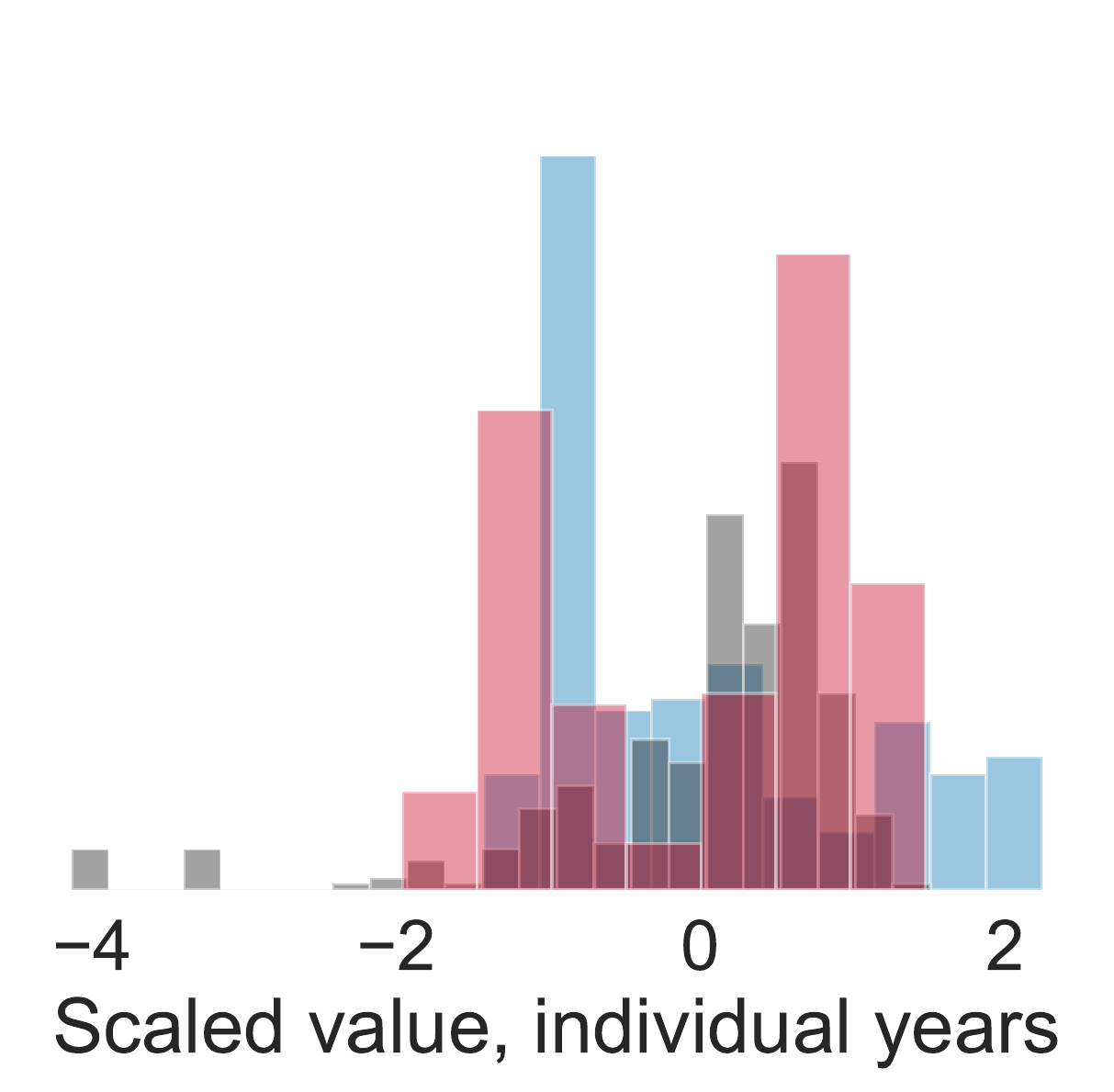}
     \caption{}\label{fig:dist3}
   \end{subfigure}
   \caption{Distribution profile for water values un-scaled (a), scaled together (b) and individually (c) } \label{fig:distribution}
\end{figure*}

\sr{Time-dependent data may have trends that makes a simple scaling meaningless, in particular when training on historic data and applying to new data.}
Fig. \ref{fig:dist1} illustrate the distribution of water values for \sr{each of} the three years of data. Fig. \ref{fig:dist2} describe the results after applying standard scaling for all years \sr{jointly}. 
\sr{The distributions differ significantly between the three years}. If two years with relative low values such as 2016 and 2017 are used for training, 
the values are not representative for the 2018 test data, \sr{which is dominated by} higher values. 
\sr{If we instead scale the water values per year, the resulting distributions become very similar, as shown in \ref{fig:dist3}}. Only scalings of water values have been illustrated in the figures, but the same pattern can be observed for several of the input variables. Scaling on individual years has therefor been \sr{applied to all variables} in this analysis. 
In real life application, \sr{it is a challenge}  that we are not able to scale the daily input variables with the yearly average since this information first is available after the year \sr{has passed}. An \sr{implementable} alternative \sr{would be to} to scale the input variable with the 365-days rolling average.

\begin{figure}[pos=htbp]
    \centering
    \includegraphics[width=0.8\columnwidth]{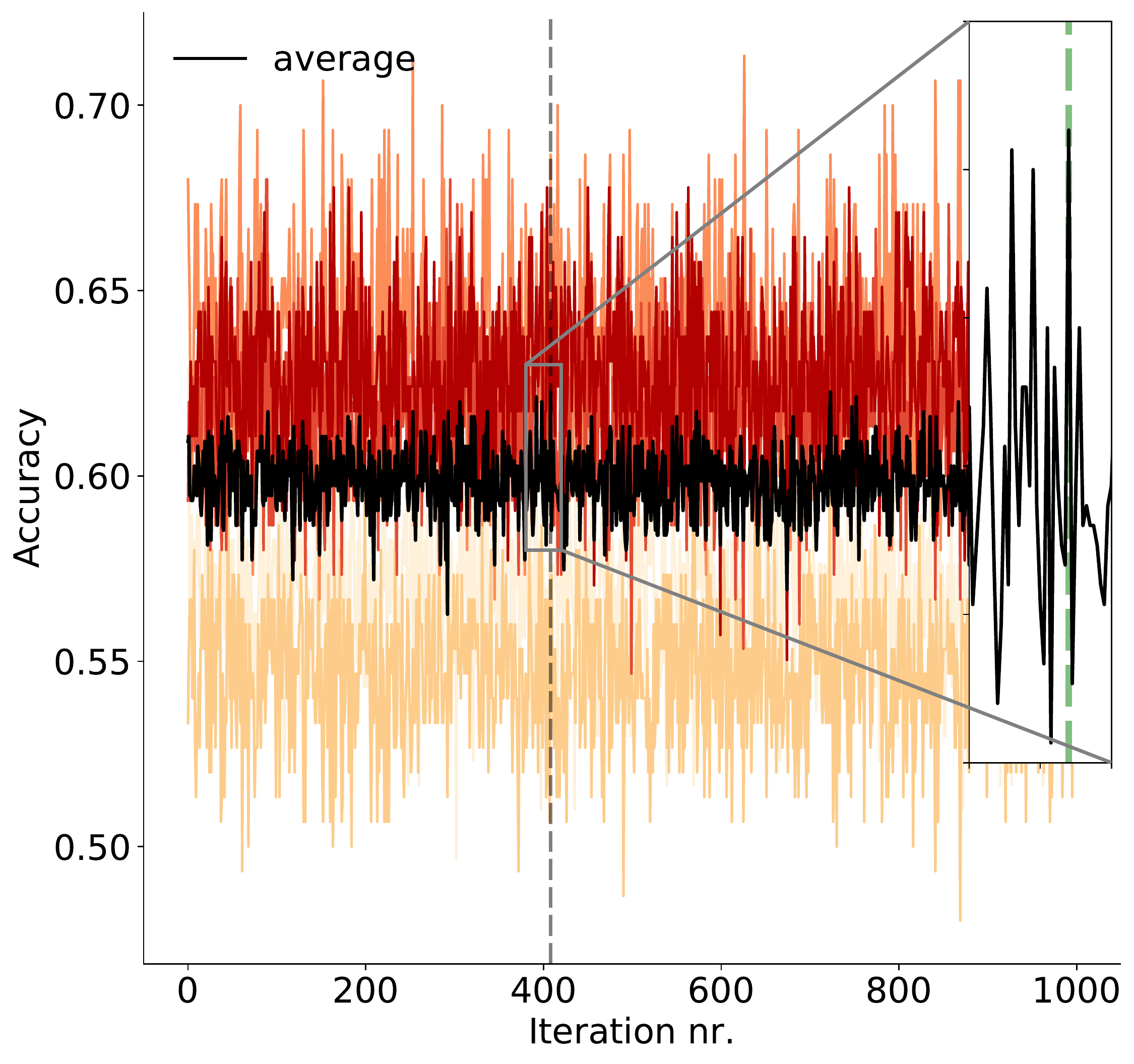}
    \caption{Hyper parameter tuning over 1000 iterations of random combination of parameters} \label{fig:hypparamtuning}
\end{figure}

\subsubsection{Hyperparameter tuning}
There are several hyperparameters in the XGBoost model that can \srs{be tuned}. 
We have applied the randomized search functionality in sklearn \cite{scikit-learn} on a selection of parameters \sr{listed in Appendix \ref{appendix-sec1} along with their final values}. 
The resulting hyperparameter values after tuning vary both depending on the number of features that are included, and also as a result of random selection of training and test data.


To avoid that the parameters are \sr{over-}optimised towards \sr{the training data}, 
\sr{the} training data can \sr{be} further split into sub-sets of data. A validation data set is a subset of the training data, 
that is used \sr{during the training to describe the evaluation of model variation due to changes in hyperparameter values and data preparation.} 


In our case, we have split the randomly selected training set into 5 sub-sets of data of approximately equal size normally defined as folds. The first fold is treated as a validation set, and the method is fit on the remaining k minus 1 folds \cite{james2014introduction}. These folds are then evaluated with 1000 set of randomly generated 
hyperparameters. The resulting accuracy obtained for each fold for every combination of hyperparameters is plotted together with the average value in Fig. \ref{fig:hypparamtuning}. \sr{Firstly, we observe a spread in accuracy within each fold illustrating the importance of tuning the hyperparameters. Further inspection of the optimization output reveals a complicated parameter landscape with multiple local minima. We also observe a spread in accuracy between each fold, indicating that the folds may not be representative for the full data set. This could be improved with additional data.  } The range in accuracy observed indicate that a fairly large spread can be expected when the model is applied on randomly generated training and test data. Finally, the hyper-parameters associated with the best average score which can be observed around iteration nr. 400 in Fig. \ref{fig:hypparamtuning}, is selected as the parameters used further in the analysis.


\subsubsection{Data dimensionality reduction} \label{sssec:DataDimensionalityReduction} 
In machine learning problems, it is tempting to include all parameters which are assumed to have an impact on the 
results. However, introducing too many features might introduce noise and distort the solution \cite{techDataDimRed}.
Reducing the number of features can be done manually by inspecting data and removing parameters with low correlation with the 
results, or high correlation with another parameter making one of them redundant. This topic is addressed in \ref{ssect:Correlation}. 

XGBoost will also by default remove variables that do not have sufficient impact on the model performance. The \sr{variables remaining after hyperparameter tuning and model optimization may be a subset of the original input variables} 
, and the result is an optimized combination of variables and parameters. When expanding the model with \sr{additional} variables, the time used for parameter tuning increases considerably. The question is if we prior to model fine-tuning, can reduce the complexity by removing variables that have limited effect on the results. \sr{Here we have used} a rule-based algorithm where \sr{model} performance is evaluated as we gradually remove features with the lowest impact on the classification results \sr{for a fixed set of hyperparameters}. \sr{The feature importance has been evaluated as the} ``gain'' computed from the decision trees in XGBoost. \sr{The gain} \cite{Gainfunc} is defined as the improvement in accuracy from adding a split on a given feature to a branch in the classification tree. Before adding the new split, there were some wrongly classified elements, but after adding the split there is now two branches which are more accurate. 

We have refit the model while subsequently removing the feature with the lowest gain. Afterwards, we select the feature combination with the highest accuracy as the best set of features. The process will be referred to as GAINS in this article.

The GAINS method has been applied on the complex data set described in Sec. \ref{sssec:complex}. Before applying GAINS, a rough screening of hyperparameters over 10 iterations is conducted. \sr{The number of input variables is then reduced based on feature importance, and the} 
final hyper-parameters are 
calculated over 1000 iterations on a smaller set of input variables.

\subsubsection{Feature importance and explanations} \label{sssec:FeatureImportance} 
Another useful tool when evaluating the relevance and importance of different features is SHapley Additive exPlantions (SHAP). \sr{The framework interprets a target model by applying Shapley game theory for how a reward given to a team should be distributed between the individual players based on their contributions \cite{SHAPpaper}.} The features are interpreted as ``contributors'', and the prediction task corresponds to the ``game''. The ``reward'' is the actual prediction minus the result from the explanation model.
The underlying idea is to take a complex model, which has learnt global non-linear patterns in the data, and break it down into lots of local linear models which describe individual data points.
\begin{figure*}[width=1.9\linewidth] 
\centering
   \begin{subfigure}{0.49\linewidth} \centering
     \includegraphics[scale=0.5]{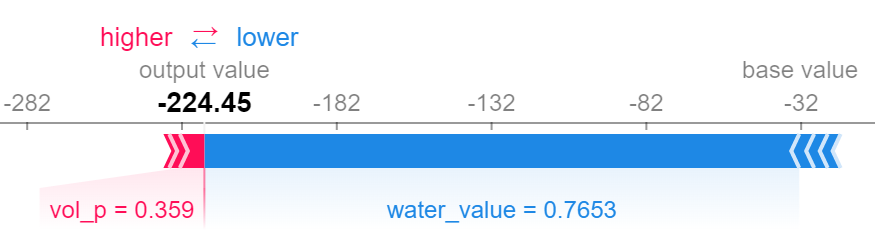}
     \caption{Individual predictions}\label{fig:SHAPfigA}
   \end{subfigure}
   \begin{subfigure}{0.49\linewidth} \centering
     \includegraphics[scale=0.5]{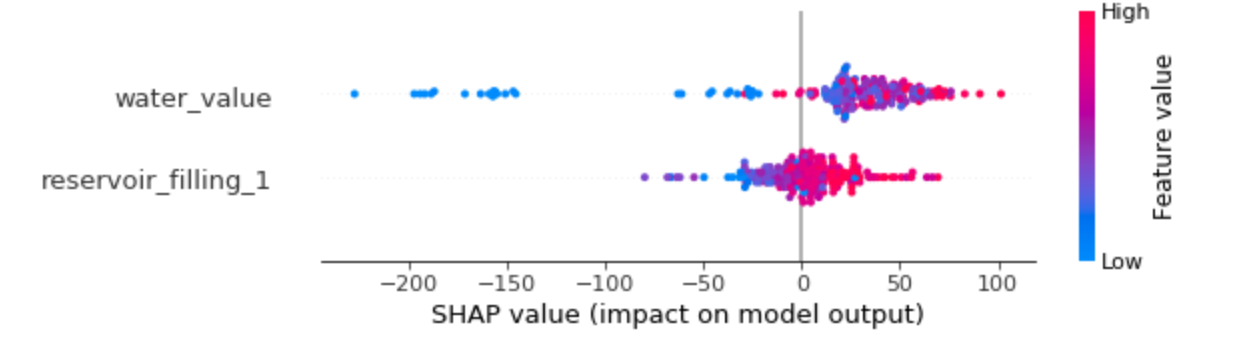}
     \caption{Aggregation of shap values for the two most important features}\label{fig:SHAPfigB}
   \end{subfigure}
   \caption{SHAP values on individual and aggregated feature values } \label{fig:twofigs}
\end{figure*}

The shap values can be interpreted either for individual predictions or for the entire sample. Fig \ref{fig:SHAPfigA} shows one sample 
classified using the regression model described in \ref{ssec:regression}. The output value is the prediction for that observation which in this case is $-224.45$ and consequently the sample is classified as stochastic. The base value of $-32$ is the value that would be predicted if we did not know any other features for the current output. This is the mean value for ``strategy gap'' for all samples in the training set. This is logical, since in lack of other information, we would predict the average value for any new samples. Interestingly, if no other information is available, we would classify the sample as stochastic. The red and blue arrows illustrate how adding other feature values push the prediction towards higher or lower values. Red indicate a push towards higher values and deterministic bidding, and blue is a push towards lower values and stochastic bidding. In this sample there is a clear push towards stochastic bidding mainly driven by the water value. The only other feature with any significant impact on this sample is the volatility for today's price which is pushing in direction of deterministic bidding.




Fig \ref{fig:SHAPfigB} illustrates an aggregation of values for two of the most important features in the regression problem. Variables are ranked in descending order of feature importance. The shap values can be interpreted as ``odds'' i.e. what is the probability of ``winning''/predict the higher/correct value which in our case what is the probability of predicting 1. Which again means predicting that the deterministic model is best. So for high (low) shap values that specific variable is contributing to increase (decrease) the probability of predicting deterministic (stochastic) bidding. This is illustrated by the horizontal location, which shows whether the effect of that value is associated with a higher or lower prediction. The color of the individual sample shows whether that variable is high (in red) or low (in blue) for that observation. If the colors are split similarly to the shap values, the relation is simple, otherwise it's probably complex and dependent on multiple variables.
A low water value has a large and negative impact on the strategy gap pushing in favour of stochastic bidding, while a high water value favours deterministic bidding. The “low” comes from the blue color, and the “negative” impact is shown on the X-axis.

The computation of shap values assume independence between the features. Our features are not independent, which may affect the ranking somewhat \cite{Aas2019ExplainingIP}. However, since shap primarily is used to qualitatively understand and illustrate potential impact of different features in this article, and  selecting features is done in combination  with other methods such as GAINS, the fact that some of the variables are dependant play a minor role in relation to prediction accuracy.



\subsection{Prediction}  \label{sec:Prediction}
\begin{figure*}[width=1.9\linewidth,pos=H]
\centering
\begin{minipage}{0.5cm}
\rotatebox{90}{\textcolor{red}{Predictions\hspace{4cm}}}
\end{minipage}%
   \begin{subfigure}{0.9\linewidth} \centering
     \includegraphics[scale=0.6]{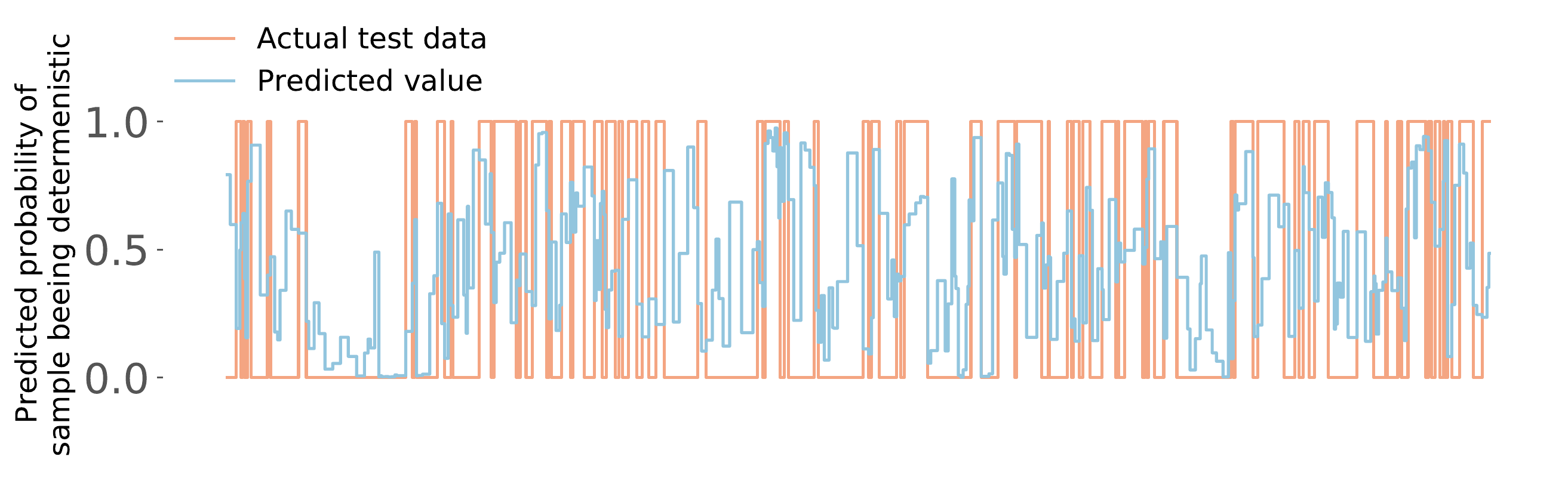}
     \caption{Classification results}\label{fig:classification}
   \end{subfigure}
   \\
   \begin{subfigure}{0.9\linewidth} \centering
     \includegraphics[scale=0.58]{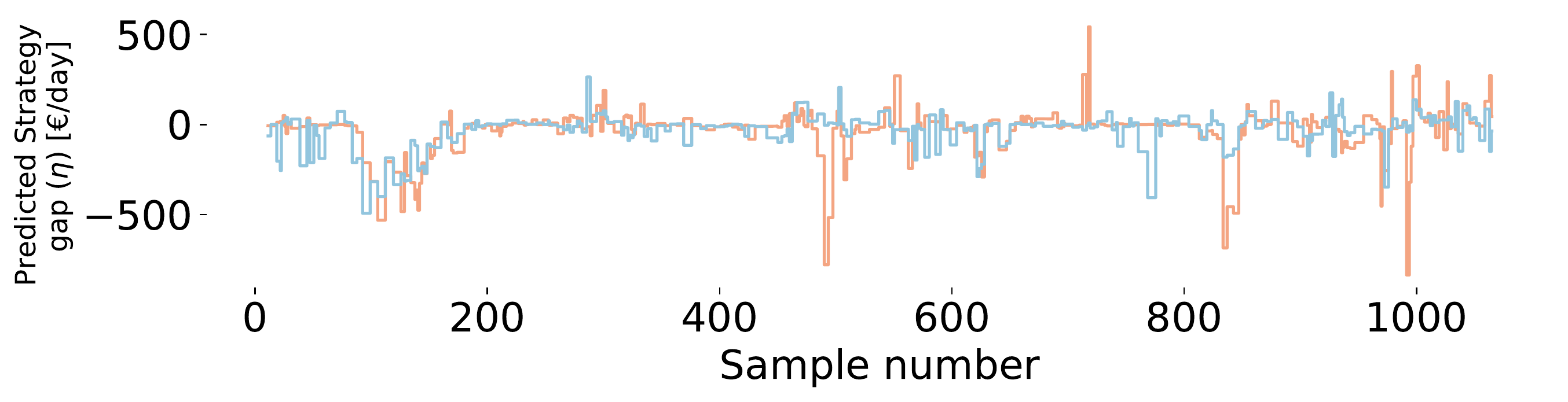}
     \caption{Regression results}\label{fig:regression}
   \end{subfigure}
   \caption{Predictions with classification and regression for 
   simple model} 
   \label{fig:classandreg}
\end{figure*}

The main target of the classification problem is to predict which bidding strategy to use for the next day. Zeros indicate that a stochastic strategy should be chosen, while ones indicate a deterministic strategy. Fig. \ref{fig:classification} illustrate predictions from the \hr{classification} model together with actual historic ``Best strategy'' from the test data. 
\sr{A nice attribute of the XGBoost classificator with binary logistic objective is that the outcome is given as probabilities and not pure classification}. This means that a figure close to zero gives a high probability of the sample being stochastic, while a figure just beneath 0.5 still classifies as stochastic, but with lower probability. 

When applying a \sr{single-output} regression model, the output 
will be the strategy gap representing the predicted difference in value for choosing one strategy in favour of another. It has previously been explained that a negative strategy gap indicates stochastic bidding, while positive numbers indicate deterministic bidding. While the classification model gives probabilities between zero and one, the regression model will span out values in a much larger range capturing the values that are at stake \sr{(after inversion of any applied scaling)}. Fig. \ref{fig:regression} illustrate the predicted strategy gap from the XGBoost regression model compared with the observations in the test data. The target is still to decide for one bidding strategy for the next day. The results from the regression model must therefore be transformed to binary recommendations \sr{of} either a stochastic or a deterministic strategy. 

Both the classification and regression model open for the possibility of applying thresholds on the results/ probabilities. As an example, when applying a threshold of 0.5 in the classification model, all values above will be classified as deterministic, while all below will be classified as stochastic. If we choose a threshold of 0.4, and 0.6, only samples with values under or over this value will be classified into a category. Values in between will not be classified. The accuracy when applying a tighter threshold will typically increase, however with the price of an increased number of unclassified days. In both these cases the model is risk-neutral when deciding on strategy. It is also possible to have a skewed threshold where we are more risk averse in relation to choosing one strategy opposed to another. If one strategy is to be selected, previous analysis \cite{Riddervold8916227} have shown that a stochastic strategy should be chosen. If a threshold of for instance 0.6 is applied, a clear tendency in direction of deterministic bidding is required before this strategy is chosen in favour of a stochastic strategy. Applying thresholds might have an even higher impact on regression models. While the probabilities associated with prediction in the classification models only are proxies for the importance of choosing the correct strategy for a selected day, the regression model indicate directly the potential losses associated with choosing the wrong strategy. A strategy focusing on classification of samples with major cost impact could be a viable solution. 



\subsubsection{Evaluation of model performance}

To monitor performance of the model, two measures are introduced:

\begin{equation}\label{eq:3}
\text{\bf{$A$}} = \frac
                {\textit{Number of correct predictions}}
                {\textit{Total number of predictions}}
\end{equation}

\begin{equation}\label{eq:4}
    {\delta}{_{realistic}} = (\overline{\beta}_{gap} - \overline{\beta}_{gap,opt}) /  \overline{\beta}_{gap,opt}
\end{equation}



Even though high classification accuracy is an important target, it is not really the main objective related to strategy selection. The main objective is to reduce the average performance gap compared to a model where the best bidding strategy is selected every day. In this sense, identifying the correct strategy on dates where there is a large performance gap between the two strategies, previously defined as the strategy gap, will be more important. This is quantified as $\overline{\beta}_{gap,opt}$, which represents the average performance gap from the optimal plan for all samples in the test data if the optimal strategy had been selected, and $\overline{\beta}_{gap}$ represents the average performance gap from the optimal plan for all classified samples. 

$\delta_{realistic}$ is then a measure of how far, we are from the optimal bidding strategy. It will be referred to as ``Realistic Performance Gap''. It is measured in percent deviation from the optimal bidding strategy. 

The model is designed to suggest an optimal bidding strategy, and performance should therefore be measured against a benchmark where an optimal bidding strategy is selected for every day, and not against an optimal plan with perfect foresight on prices and inflow. This is why the measure of $\delta_{realistic}$ is introduced, rather than using $\overline{\beta}_{gap}$ as performance measure. 

Assuming that the consequences of wrongly classifying samples are normally distributed between stochastic and deterministic strategies, and that there is a fairly equal split between when the two strategies perform best, the accuracy represents a good measure for model performance. In other cases, for instance in medical  diagnostics, the consequence of failing to identify disease of a sick person might be much higher than the cost of sending a healthy person to more tests. In this case we have to use more sophisticated methods for evaluating performance of the model \cite{Confmatrix}.

\begin{figure}[pos=htbp]
    \centering
    \includegraphics[width=0.8\columnwidth]{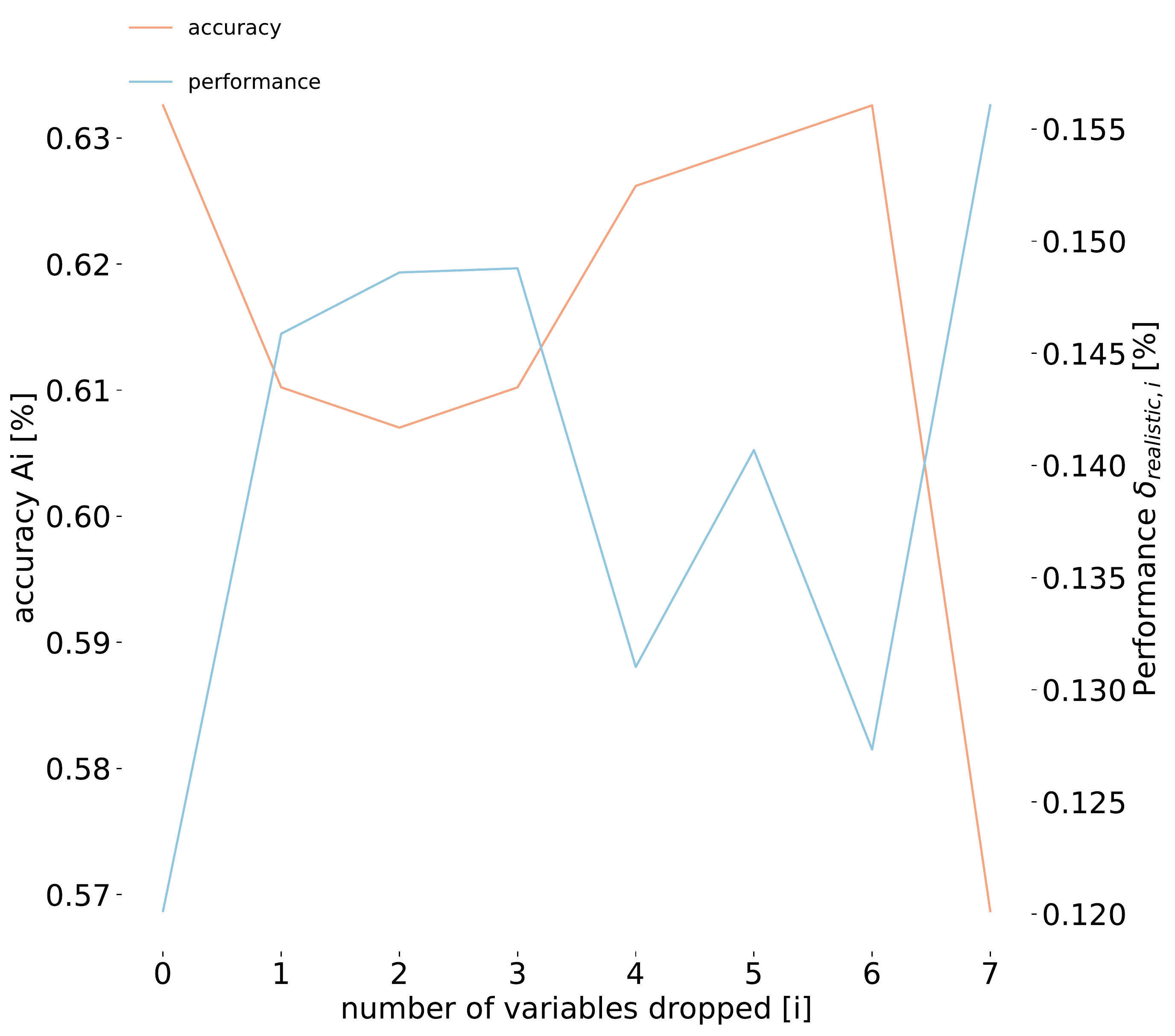}
    \caption{Variable reduction} \label{fig:varreduction}
\end{figure}

In fig. \ref{fig:varreduction} the two measures of \textit{A} and $\delta_{realistic,i}$ are plotted as a function of number of variables that are removed in the classification problem. In this example, the best accuracy is obtained with the eight original variables. 
It can also be observed that this coincides 
with when the gap to a optimal bidding model is the lowest. 

\section{Concept for application and case study} \label{casestudy}
An important prerequisite when building a learning model, is the access to data, and specifically in this case, historic bidding performance. To benefit from previous work, this analysis is based on the results obtained in \cite{Riddervold8916227} where bidding performance for a river system located in South-Western Norway has been analysed for 2018. Given the large variations that are associated with operation of a hydro-based system, one year of data is considered to be too short. The period evaluated has therefore been extend with additional simulations for the years 2016 and 2017.

The framework for how a boosting algorithm can be used to indicate the optimal bidding strategy for any given day has been presented in section \ref{ML}.
The examples used to illustrate the learning process are all linked to the real life case discussed in this section. The data used for illustration purposes in the first section are based on the eight original variables. When the amount of variables increase considerably, so does the complexity in visualisation and graphical interpretation. Introducing hourly variables on prices, price prognosis, and bid-ask sensitivities increases the number of variables to over 100. This is still considered to be a small problem in machine learning terms (e.g. genetic problems deal with several hundred thousand features \cite{XGboostgenes}), but correlation matrices, tree-structures, feature-importance diagrams etc. will contain more information than can be practically visualized. 
The results from the processes can rather be summarized with presentation of the resulting parameters after tuning and training, and the results of the classification on the test data.

Table \ref{tab:resultsall} summarizes the results from different modelling approaches on various splits of the available data.

\begin{table*}[width=1.9\linewidth,pos=h] 
   \centering
   \caption{Results for predicted accuracy (A) and realistic performance gap ($\delta_{realistic}$) with different modelling approaches}
   \label{tab:resultsall}
   \begin{tabular}{p{1.7cm} p{6.5cm} p{0.9cm} p{1.1cm} p{1cm} p{2cm} } 
     \hline
    Model & description & A & $\delta_{realistic}$ & $A\textsuperscript{*}_{mean}$ & $A\textsuperscript{*}_{std}$  \\ \hline
    (1) XSCR & Xgboost Simple Classification Random & 0.62 & 0.18 & 0.61 & 0.03 \\
    (2) XSCS & Xgboost Simple Classification Sequential & 0.56 & 0.21 & 0.55 & 0.03 \\ 
    (3) XCCR & Xgboost Complex Classification Random & 0.63 & 0.15 & 0.63 & 0.03 \\
    (4) XSRR & Xgboost Simple Regression Random & 0.60 & 0.15 & 0.60 & 0.03 \\
    (5) NNRR & Neural Network Regression Random & 0.59 & 0.19 & NA & NA \\
    (6) RNRS & Recurrent Neural Network Regression Sequential & 0.57 & 0.22 & NA  & NA \\

     & & &   \\
\end{tabular}
\end{table*}

\subsection{Classification}
\label{ssection:Classification}
\subsubsection{Simple model - Original variables}
Two main approaches \sr{for selecting training and test samples} have been investigated when applying the classification model, sequential and random split. 

\sr{We observe a variation in performance for different random seeds. This indicates noise in the data, and that finding a }
universal set of variables and hyperparameters that fit well with all data, might be difficult. The result for the random sampling in Table \ref{tab:resultsall} is therefor when the model is fitted using the eight original variables. In order to determine the variance of possible outputs on the random sample, we bootstrap \cite{Statistics_out-of-bagestimation} the test sample and determine the performance on each sub-sample individually. The mean value \hr{($A\textsuperscript{*}_{mean}$)} and standard deviation\hr{($A\textsuperscript{*}_{std}$ )} after boostrapping on 100 
samples is also indicated in the same table.



As expected, applying the model on sequential data give poorer \sr{accuracy performance than for random samples}. When 2016 and 2017 data are used to predict 2018 strategies, the results are barely better than random guessing. When investigation of 2016 and 2017 data in addition reveal that stochastic bidding is best 58\% of the time, and the amount of stochastic samples in 2018 data is 52\%,  purely applying a binomial selection with a probability factor for stochastic results equal to 0.58 would increase the accuracy to above 50\% without taking into account any variables.

One would however, expect that the results from sequential splitting would approach the accuracy obtained by random split as more data will become available. 

The random sampling 
give a\sr{n accuracy} around 62\% , but also here the training and test set consists of an average of 54\% stochastic samples, so we are only increasing accuracy by 8\% compared to applying a binomial selection without any variable input.  

An interesting observation is that even though there is an average of 18\% realistic performance gap ($\delta$) for the prediction model when applying random sampling, the gap is still 3\% better than when applying a pure stochastic strategy for all days in the evaluation period. The benefit in addition to a marginal better performance, is the reduced computational capacity required to perform the bidding process, since we reduce the number of days when the stochastic bidding process must be used by almost 40\%.
Relying solely on a deterministic bidding strategy is not a good idea, since this would give a realistic performance gap ($\delta$) of 43\%.







\subsubsection{Complex model -  Variable extension and reduction} \label{sssec:complex}
The results from the boosting model when applying the limited set of original variables will potentially give the power producer better insight in relation to choosing the correct strategy, but the accuracy of the model still seems to be low compared to what is required from a good prediction model. \sr{Will} increasing the amount of variables and complexity in the model increase the performance? 
In addition to daily average values associated with price and volatility, all hourly prices for both issue date and prognosis for value date are included. The hypothesis is that the classification model interprets the inherent volatility and importance of individual hours when predicting the bidding strategy. Bid-ask sensitivities as described in \ref{sssec:Bidaskcurves} have also been added. Finally, additional features such as month, year, day, performance of similar week-days, performance gap for day minus one, rate of change for reservoir filling and difference between price and water value are included in the model. The extended model now consist of 113 variable compared to the 8 original variables. A correlation matrix and intuition would reveal that several of the new and old variables are correlated and could potentially be removed. 

\sr{We apply the GAIN-loop and SHAP analysis to the same data sample as used for the eight original variables. Due to the computational requirements from increasing the number of features, the initial hyperparameter tuning is reduced to 10 iterations. This is both to reduce the time used in the tuning process, but also to avoid that the parameters are customized to the data-set consisting of all variables.
After the parameter tuning process, only 15 of the original 113 variables are still present in the model, and the GAIN-loop further reduces the number to 11.}



\begin{figure}[pos=htbp] 
\centering
     \includegraphics[scale=0.26]{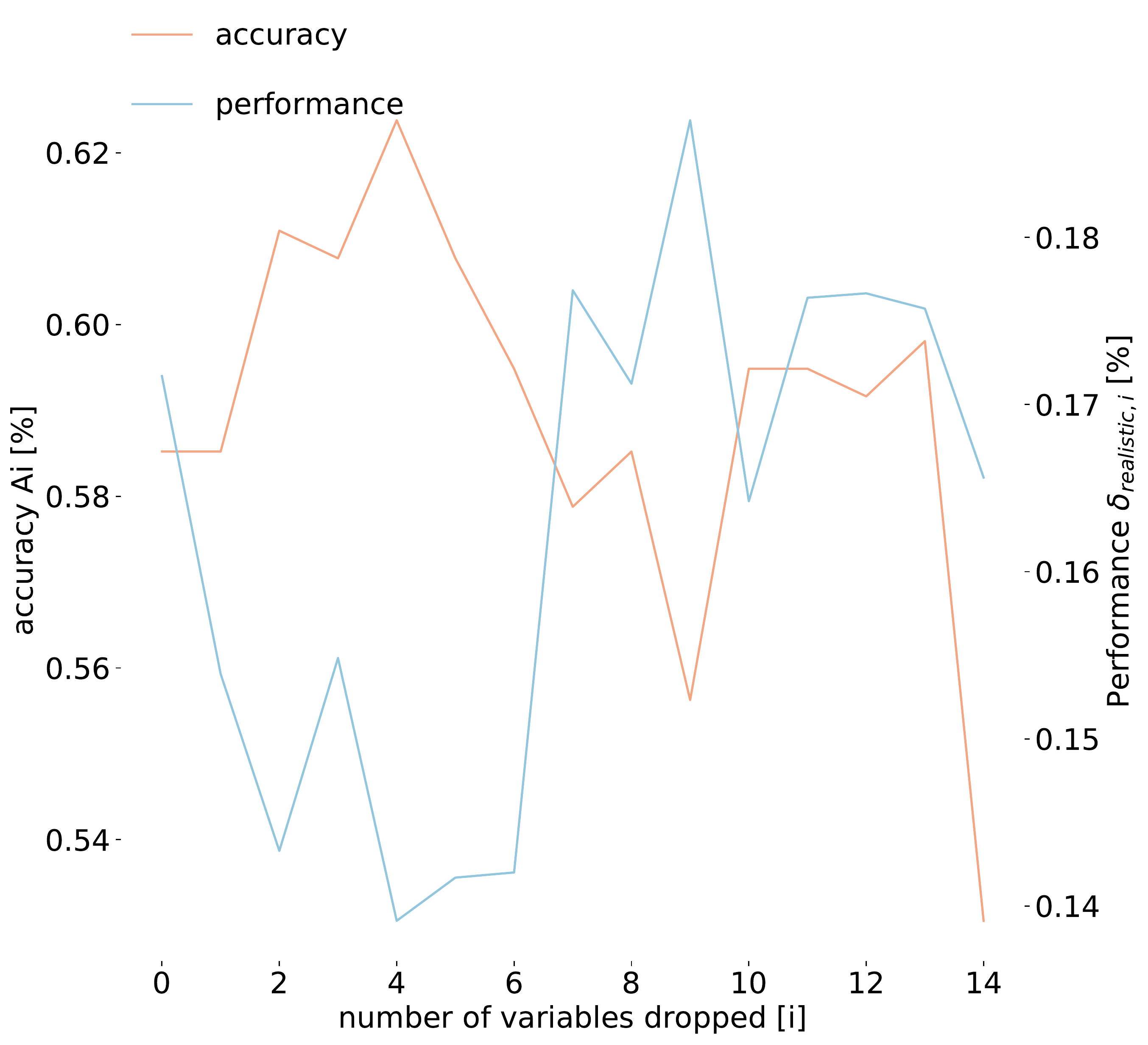}
     \caption{Variable reduction for complex model}\label{fig:gainresults}
     
     \end{figure}

Fig. \ref{fig:gainresults} illustrate how the accuracy gradually increase as more variables are removed until we reach the highest accuracy at \textit{i} = 11. The realistic performance gap($\delta$) is also low at this point.

\subsubsection{Comparison of GAINS and SHAP}

\begin{figure}[pos=htbp] \centering
     \includegraphics[scale=0.28]{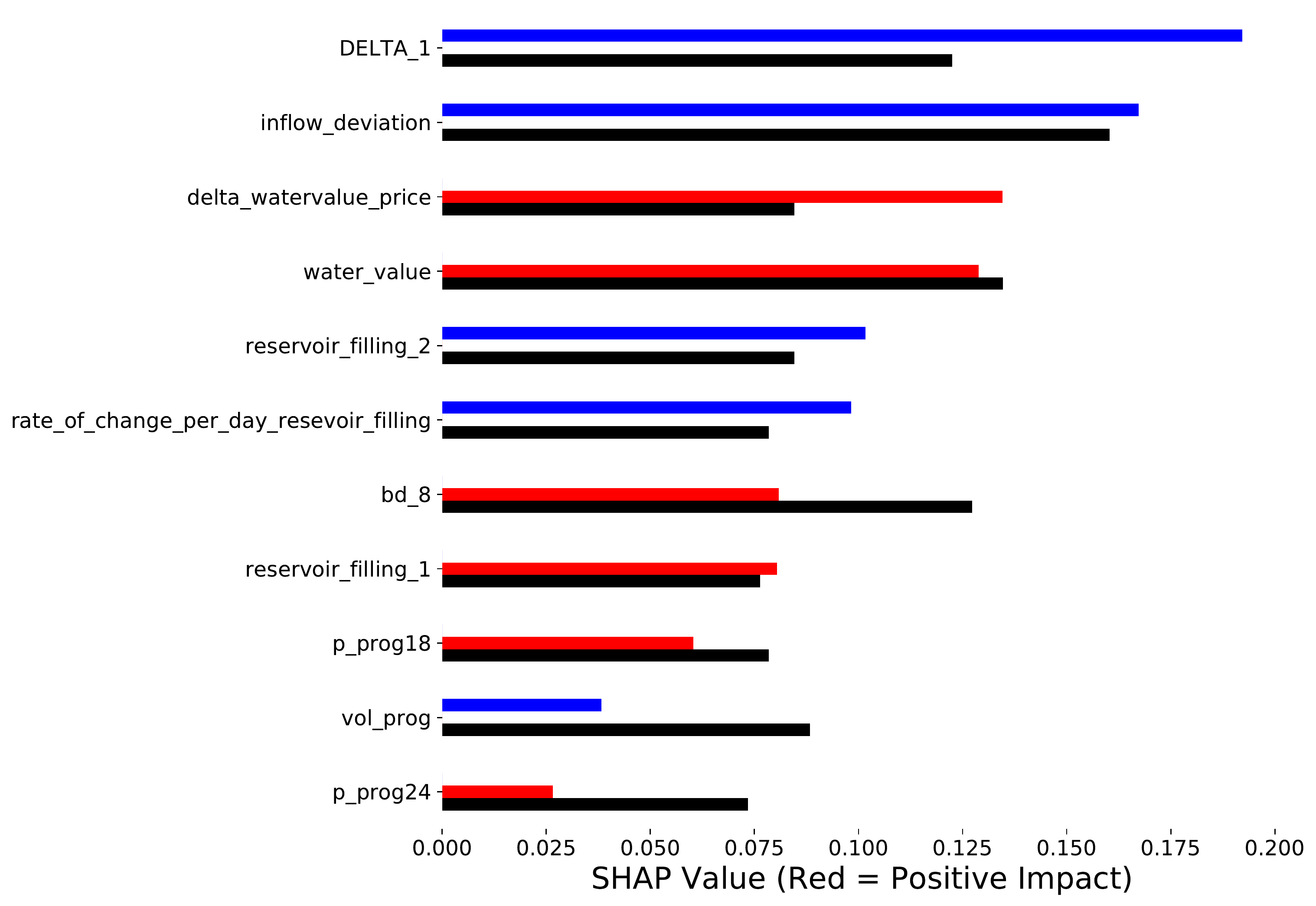}
     \caption{Top 15 SHAP values sorted compared with features importance from GAINS \sr{Should we explain the variables here?}}\label{fig:SHAPvaluessorted}
\end{figure}

\sr{SHAP and GAINS values carry different information \cite{SHAPpaper}. In Fig.\ref{fig:SHAPvaluessorted}, the feature importance from GAINS is displayed (in black) together with the resulting values from the SHAP analysis (red/blue for positive/negative impact). Feature importance from the GAIN-loop have been normalised to the average value of the SHAP samples to facilitate for better comparison.}


The SHAP values are sorted with the most important variables at the top. Red color means a feature is positively correlated with the target variable, while blue color indicate negative correlation. High water values 
\sr{correlates with} deterministic bidding. 
Inflow deviation shown in blue indicate that low values for this variable favours deterministic bidding. The same is the case for strategy gap for day minus 1 (DELTA\_1). 
We already know that a high value for strategy gap indicates deterministic bidding. The SHAP value further implies that if deterministic bidding was the best for day minus one (issue date), there is an increased probability for deterministic bidding also is the best strategy for the day we investigate bidding for (value date). 

There is clearly a deviation between the two methods. SHAP seems to emphasize the importance of the strategy gap observed for day minus one to a larger extent than feature importance from GAIN. On the other hand, GAINS give higher relative scores for hourly variables such as bid sensitivity in hour 8 (bd\_8) and price prognosis for hour 24 (p\_prog24). 

When the SHAP analysis is \sr{performed} after a considerable number of variables are removed, which is the case 
in Fig. \ref{fig:SHAPvaluessorted}, we have lost information about the potential effect of the \sr{filtered} variables. One would  assume that these anyhow contribute with limited information, but a test applying SHAP on an un-reduced model, indicate that two variables, (vol\_roll\_2) and (bd\_1), also appear among the top eight variables in the SHAP analysis. These are therefor also included when moving forward to the next step in the analysis.
  
Water value, inflow deviation and reservoir filling 2 which were part of the eigth original variables are listed as important in both SHAP and GAINS. In addition, features such as strategy gap for value day minus one, delta water value and price, rate of change for reservoir filling seem to have significant impact on the resulting predictions.
Furthermore, GAINS seem to emphasize the importance of the bid curves for hour 8 (bd\textunderscore8) and the volatility in the next days prognosis (vol\_prog). 
An interesting observation is that from the 96 hourly variables, we are now left with two. These are most likely chosen by the model because they represent the possibility for low prices during night (bd\_1) and morning peak prices (bd\_8).      

Exactly where the cut-off should be in relation to how many variables should be brought further to the next step in the analysis is difficult to measure, but three additional variables (reservoir\_filling\_1, p\_prog18, and p\_prog24) listed in Fig.\ref{fig:SHAPvaluessorted} are removed manually since these neither seem to be significant in GAINS nor SHAP. 

\hr{The 10 final variables are the 8 remaining variables in Fig.\ref{fig:SHAPvaluessorted}, in addition to (vol\_roll\_2) and (bd\_1).}
The next step is to re-run the parameter tuning process with the 10 variables over 1000 iterations, and test the model on \sr{randomly sampled training and test data}. The same random seed as used for evaluating the original model with eight variables is used to illustrate the performance on an identical set of data.

The accuracy using these 10 variables turn out to be 63\% which is ~ 1\% better than the result from the simple model where the eight original variables are used. Given that the 113 new variables are reduced to 10, it is a clear indication that most of the new variables contain little or no additional information relevant for improving the prediction accuracy. In contrary, including new variables might just as well introduce noise to the problem and distort the solution as explained in section \ref{sssec:DataDimensionalityReduction}.

Two main conclusions can be drawn from the comparison between SHAP and GAINS, and the re-run of the model with reduced number of variables. 
Firstly, variable reduction based on GAIN \hr {and SHAP} might not necessary end up giving us the combination of variables that give the best prediction. This is because the potential interaction between variables are not sufficiently addressed in this process \cite{shapexample} 
Secondly, there is substantial amount of noise in the data \sr{leading to high sensitivity of the final results to the choice of variables and the composition of traning and test data.} 
This leaves us with a relatively weak prediction model. \sr{However, the situation can be improved with more/less noisy data.}

\subsection{Regression} \label{ssec:regression}

The classification process described in \ref{ssection:Classification} is primarily designed and optimized towards obtaining the best prediction accuracy for stochastic and deterministic bidding. The best accuracy correlate (albeit weakly) with when the lowest realistic performance gap ($\delta$) is observed, but it can be seen that there are instances where a low ($\delta$) could have been obtained at the expense of poorer overall accuracy.
In Fig.\ref{fig:gainresults} the performance gap ($\delta$) barely change when the number of variables dropped increase from four to six, but the accuracy is reduced by almost 4\% from 62\% to 58\%

If the primary target is to obtain as low ($\delta$) as possible, several adjustments to the model could be made. 
\sr{First of all, the classification model has been trained with a binary logistic objective, which specifically optimizes the accuracy regardless of the resulting gap. Fitting for the actual gap values with e.g. a ``Mean Squared Error'' (MSE) objective would penalise large gaps more than small gaps.}

It is also possible to introduce more categories and perform a multi-class classification. One could for instance use five groups where they represent ranges for strategy gap. One group could be samples with strategy-gap's in the range -400 to -200. These are samples which we obviously would like to predict as stochastic, and additional weight on this group could be incorporated in the objective (loss) function. 


The input data has already been classified in stochastic and deterministic based on the positive/negative values of the strategy gap. 
To take into account the importance of the numeric values of the strategy gap, an alternative method is to perform regression on the strategy gap directly, and rather classify the best strategy with simple heuristics after predicting the strategy gap.
The reason for not pursing this method as our primary approach, is that the relative limited amount of data, and amount of noise in the input data, which lead to a model with relative poor performance on predicting the strategy gap. However, if we assume that we are only interested in if the model predicts a number that is higher or lower than zero, the regression model might actually be able to incorporate an important aspect not sufficiently addressed in the classification. \sr{With a mean squared error objective, the GAIN-loop can focus on variables that contribute to reducing the gap}. The criteria for selecting the best set of variables can therefore be based on the strategy gap, rather than the classification accuracy. MSE is given by eq.\ref{eq:4}, where $\eta_i$ is the strategy gap for sample i and $\hat\eta_i$ is the predicted value.
\hr{
\begin{equation}\label{eq:5}
MSE = \dfrac{1}{n}\sum\limits_{t=1}^{n}(\eta_i - \hat{\eta}_i)^2
\end{equation}
}
Given the relatively small difference in results obtained by introducing the increased complexity in the complex model, regression results have only been tested on the simple model, and the results can be seen in table \ref{tab:resultsall}. As expected, the classification accuracy is 2\% worse than for the model optimised towards classification, but the performance gap for the model is actually 3\% better. 





\subsubsection{Neural Networks}
\label{sssec:neuralnetwork}
The XGBoost algorithm used in the analysis in the previous sections was chosen for two main reasons. First it is a method allowing for transparency in relation to tracking importance of the different variables in the prediction process, and thereby avoiding the ``black-box'' perception that is associated with neural networks. The second argument is that the algorithm is suitable for both classification and regression. 

\begin{figure*}[pos=htbp]
\centering

\begin{subfigure}[t]{.5\textwidth}
  \centering
  \includegraphics[scale=0.45]{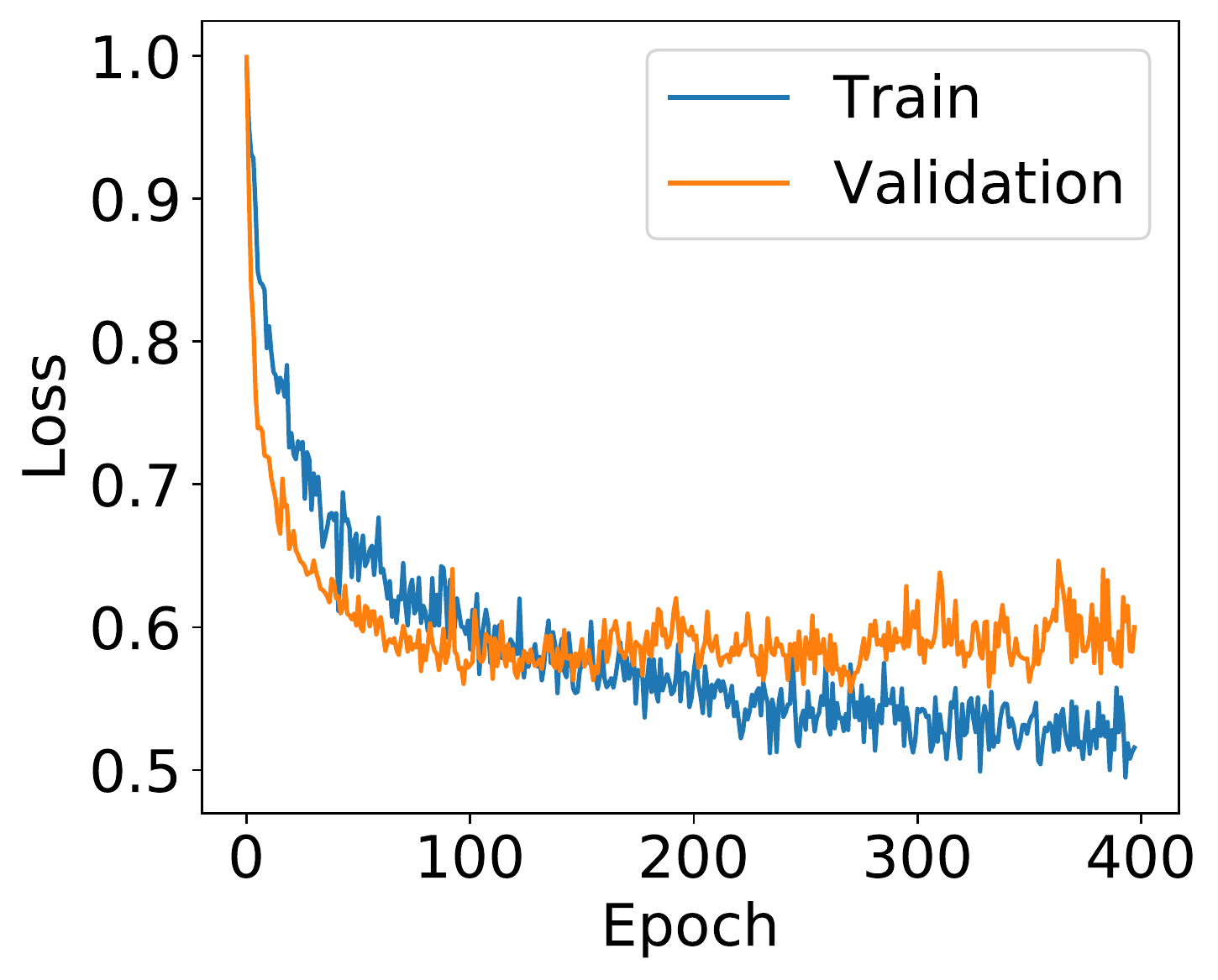}
  \caption{Well fitted Neural Network with eight variables }
  \label{fig:NNvalidationgoodfit}
\end{subfigure}%
\begin{subfigure}[t]{.5\textwidth}
  \centering
  \includegraphics[scale=0.45]{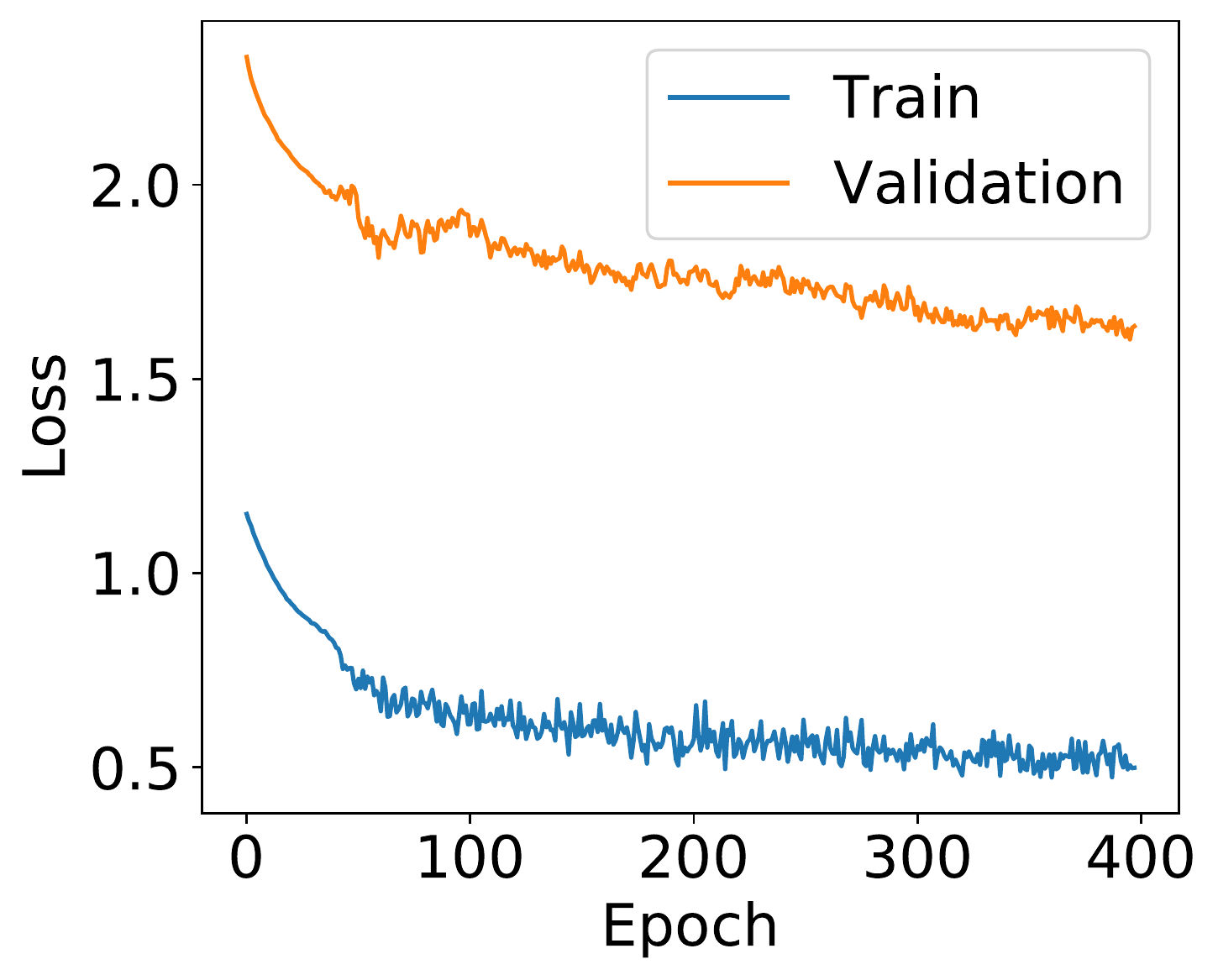}
  \caption{Poorly fitted Reccurent Neural Network}
  \label{fig:NNvalidationbadfit}
\end{subfigure}
\caption{Plot of model loss for Neural Networks }
\end{figure*}

As a benchmark for results obtained by the XGBoost regression model, a multilayer perceptrons (MLPs) model, often referred to as a classical type of neural network, has been tested on the same random sample analysed using the XGBoost model on the original set of variables (XSRR) . 
Typically, with neural networks, we seek to minimize the error the prediction model makes, \sr{also referred to as the loss.} 
When training a neural network, it is important to avoid over-fitting. This can be monitored by comparing the loss function for the training and validation sets. If the two values converge with training, the result from the model is assumed to fit well. An example of a good fit can be observed for the 150 first iterations in  Fig. \ref{fig:NNvalidationgoodfit}. This is the results from 400 iterations 
with the eight original variables. When loss for the training and validation data diverge as the number of iterations increase, the model is over-fitted. This seems to be the case if the numbers of iterations are increased beyond 150. The result in this case, is a model which is really good at representing the training data, but not necessarily good at predicting results for unseen data. The neural network benchmark is therefore established using results after 150 iterations, and the results can be found in table \ref{tab:resultsall}.
The prediction accuracy and realistic performance gap is ($\delta$) somewhat worse than for the regression using XGboost.
\hr{The Neural Network models have been manually tuned, and the resulting model architectures and hyper-parameters are listed in Appendix \ref{appendix-sec1}}

\sr{In a traditional neural network we assume that all inputs (and outputs) are independent of each other. To exploit any time dependence in the data, we have applied a recurrent neural network (RNN) on the sequential split of training and test data.  Since time-dimension is a key aspect of RNNs, it does not make sense to apply it on random selection of data.}

The most relevant test of the RNN in this case study, is to train and validate the model on the historic data in sequential order, and test in on samples that follow sequentially after the period of training. To be able to compare with a similar period investigated with the classification model, the RNN has been used with training in 2016 and 2017, and tested on 2018 data. The network has been \sr{trained over 400 epochs}.   
The results can be seen in table \ref{tab:resultsall}. \sr{Similarly to the classification results} for the same period, the accuracy is poor with only 57\% hit-rate, and a realistic performance gap ($\eta$) of 22\%. Using this model to select bidding strategy would actually give higher losses than purely choosing a static strategy with stochastic bidding. The loss for training and validation as shown in Fig.\ref{fig:NNvalidationbadfit} indicate that we have a poorly fitted model. 
\subsection{Alternative and combined approaches}
Both XGBoost and neural networks have strength and weakness that can influence the results of the analysis. A nice attribute of the XGBoost algorithm is the ability to track the importance of different features, and thereby remove features that have little or no impact on the results. A useful attribute in the neural network is the possibility predict more than one value, and combine this with a custom made loss function. A hybrid model where XGBoost is used to select the most important features, and a Neural Network is applied with the recommended set of features and a custom made loss function (CL), could be a possibility to combine the strengths of the two models. This approach has been investigated 
, and indicate that results might improve marginally for the neural network model, but the results are still poorer than for only using the XGBoost model. The custom loss function used is designed to focus on the days where the highest performance gaps can be observed. 

\hr{
\begin{equation}\label{eq:6}
CL = \dfrac{1}{n}\sum\limits_{t=1}^{n}|min(\beta_{det}, \beta_{stoch})- min(\hat\beta_{det}, \hat\beta_{stoch})|^{\it{n}}
\end{equation}
}

In comparison with Eq.\ref{eq:5}, Eq.\ref{eq:6} performs regression on the performance gap ($\Pi_s$) for each bidding strategy directly, and not on the strategy gap ($\eta_s$). This makes it possible to ``punish'' the days where one strategy perform poorly compared to another. The delta can be raised to a power of \textit{n} to emphasize the difference even more.

Antoher concept that could be an alternative to the investigated approaches is Bayesian model averaging, where multiple models can contribute to the decision support, with weights based on how sure each model is.




\section{Conclusion} \label{Conclusion}
In this article, various techniques for classification and regression have been tested on historical data representing bidding performance for a reservoir-based river-system in the Nord Pool market. The primary objective has been to investigate the possibility of predicting prior to day-ahead bidding whether stochastic or deterministic bidding would be the preferred strategy under the prevailing market and hydrological conditions. 

A simple plot of inflow deviation together with performance gap for the year 2018 as seen in Fig \ref{fig:performacegap_vs_inflow} gave an indication that periods with high inflow to a certain degree were correlated with periods with negative stategy-gaps favouring stochastic bidding. This was partly the motivation for conducting the analysis, and also one of the hypotheses in the initial selection process for variables to be used in the prediction model. 
Inflow deviation together with water value has also proved to be the variables that play the largest role in the prediction model. This has been demonstrated using two different techniques for investigating feature importance that are associated with the gradient boosting in XGboost. 

If historical data are assumed to represent future instances with sufficient accuracy, applying a prediction model trained on available data will be able to predict the correct strategy with an accuracy of 62-63\%. If the benchmark is a strategy where only stochastic bidding is performed, the prediction model will outperform this strategy, and the number of days that must be analysed with the computationally demanding stochastic model is reduced by almost 40\% .
A considerable increase in variables beyond what was assumed by to be the most important ones by domain experts, do not not seem to improve prediction accuracy considerably. A prediction model with only 62\% accuracy and a standard deviation of 3\% indicate that there could be considerable amount of noise in the data. Given that only one combination of hyper-parameters and variables can be selected in an operational environment, the risk of over-fitting a model increase when the model is fit to a large set of variables. This might favour choosing a simple model with a limited set of variables.

Neural networks have been tested to benchmark result from XGBoost, but results indicate that the initial assumption that gradient boosting models are suited for this type of analyses was well founded. Using regression rather that classification has proved to be able to reduce the performances losses, but at the expense of somewhat lower accuracy score. 

One challenge associated with the use of the prediction model arise when the model is used to analyse data that have a distribution profile on key variables that vary considerably from what is represented in the test data. This is clearly the case when the model is used to analyse a year without any samples represented in the training set. Increased amount of historic data will most likely increase the probability that training data are representative, but the volatile nature of power markets and weather does not guarantee that that this is going to be the case. 


If a power producer decides to implement a strategy selection model, including new observations and re-calibrating the model on a daily basis could be the most robust solution. 
\srs{One difference between Neural networks \hr{and} boosting methods is that Neural network can be updated on the fly, whereas boosting methods must be fully retrained when the training data changes.}
Tracking performance of a strategy selection model will be important, and the power producers should consider running the model in parallel with existing systems for a period of at least one year to account for seasonal effects. If performance of the model turns out to be poor, it is an indication that the historic information available might be insufficient to give results with sufficient quality. 



\newpage
\appendix
\section{Appendix 1}
\label{appendix-sec1}

\begin{table}[pos=hb] 
   \centering
   \caption{Model architecture and hyper-parameters}
   \label{tab:Modelarchitecture}
   \begin{tabular}{p{3cm} p{1.5cm} p{1.5cm} p{1.5cm} p{1.5cm} p{1.5cm} p{1.5cm} } 
     \hline
    Param & XSCR & XSCS & XCCR & XSRR & NNRR & RNRS\\ \hline
    Algorithm & Xgb & Xgb & Xgb & Xgb & Keras & Keras \\
     &  &  &  &  & Sequential & Sequential \\
     &  &  &  &  & Dense & LSTM \\
    objective & binary:  & binary: & binary: & reg: & mse & mse\\ 
             & logistic  & logistic & logistic & squarederror &  & \\
    learning\_rate & 0.075 & 0.075 & 0.032 & 0.092 &  & \\
    max\_depth & 4 & 4 & 4 & 3 &  &\\
    n\_estimators & 178 & 178 & 323 & 259 &  &\\
    number\_boost\_round & 9282 & 9282 & 2761 & 2509 &  &\\

    gamma & 4.26 & 4.26 & 9.68 & 4.34 & &  \\
    subsample & 0.77 & 0.77 & 0.85 & 0.64 &  & \\ 
    dropout\_frac &  &  &  &  & 0.2 & 0.3 \\
    number\_neurons1 &  &  &  &  & 50 & \\
    number\_neurons2 &  &  &  &  & 20 & \\
    lookback &  &  &  &  &  & 30 \\
    L2 &  &  &  &  & 0.005  & 0.005 \\
    kernel\_regularizers &  &  &  &  & l1(L2)  & l1(L2) \\
    activation &  &  &  &  & linear & linear \\
    optimizer &  &  &  &  & adam  & adam \\

 \hline
     & & &   \\
     

\end{tabular}
\end{table}

\newpage
\hrblank{whiteblanks}
\newpage

\bibliographystyle{cas-model2-names}


\bibliography{cas-refs}





\end{document}